\documentstyle[11pt,epsf]{article} 

\newtheorem{sect}{}[section]

\begin{document}

\title{Geometric Interpretations of Quandle 
Homology}

\author{
J. Scott Carter \\
University of South Alabama \\
Mobile, AL 36688 \\ carter@mathstat.usouthal.edu \and
Seiichi Kamada \\
Osaka City University \\
Osaka 558-8585, JAPAN\\ kamada@sci.osaka-cu.ac.jp \\
skamada@mathstat.usouthal.edu
\and
Masahico Saito \\
University of South Florida \\
Tampa, FL 33620 \\ saito@math.usf.edu
}
\maketitle
\begin{abstract}

Geometric representations of cycles in quandle homology theory are 
given in terms of colored knot diagrams. 
Abstract knot diagrams are generalized to diagrams 
with  
exceptional points
which, when colored, 
correspond to 
degenerate 
cycles. 
Bounding chains are realized, and used to obtain equivalence moves 
for homologous cycles.
The methods are applied to prove that boundary homomorphisms in 
a homology exact sequence vanish. 
\end{abstract}

\section{Introduction}

A quandle is a set with a 
self-distributive binary operation (defined below)
whose definition was motivated from knot theory. 
A (co)homology theory was defined in \cite{CJKLS} for quandles,
which is a modification of rack (co)homology defined in \cite{FRS2}. 
State-sum invariants using quandle cocycles as weights are 
defined \cite{CJKLS} and computed for important families
of classical knots and knotted surfaces \cite{CJKS1}.
Quandle homomorphisms and virtual knots are applied to this 
homology theory \cite{CJKS2}. The invariants were applied to study 
knots, 
for example, in detecting non-invertible
knotted surfaces \cite{CJKLS}.

On the other hand, knot diagrams colored by quandles can be used 
to study quandle homology groups. This view point was developed
 in \cite{FRS2,Flower,Greene}  
for rack homology and homotopy. 
In this paper,
we generalize this approach  to quandle homology theory. 
For this purpose, we generalize abstract knot diagrams \cite{KK},
and represent cycles and boundaries 
geometrically by colored abstract knot diagrams. 
Equivalence moves are given in terms of colored knot diagrams.

The paper is organized as follows. 
In Section~\ref{def} 
the  
definition 
of quandle homology 
is 
recalled.
Abstract knot diagrams 
are generalized in Section~\ref{diagsec}.
Colorings of generalized diagrams are defined (Section~\ref{colorsec}),
and used to represent cycles of quandle homology (Section~\ref{repsec}).
Equivalence moves for colored diagrams are given in Section~\ref{equivsec}.
In Section~\ref{fundsec}, the fundamental quandles for abstract knots are
 discussed. 
The boundary homomorphisms in an exact 
sequence are shown to be trivial (Section~\ref{bdryhomsec}).
Examples appear throughout the paper.

\section{Definitions of Quandle (Co)Homology} 
\label{def}

\begin{sect}{\bf Definition.\/}
{\rm 
A {\it quandle\/}, $X$, is a set with a binary operation
$\ast$ such that 

\noindent
(I. {\sc idempotency}) for any $a \in X$, $a \ast a=a$,

\noindent
(II. {\sc right-invertibility}) for any $a, b \in X$,
 there is a unique $c \in X$ such that $ a =
c\ast b$,
and 

\noindent 
(III. {\sc self-distributivity}) for any $a, b, c \in X$, we have
 $(a \ast b) \ast c =
(a \ast c) \ast (b \ast c)$.

A {\it rack\/} is a set with a binary operation that satisfies (II) and
(III).
Racks and quandles have been studied in, for example, 
\cite{Brieskorn},\cite{FR},\cite{Joyce},\cite{K&P}, and \cite{Matveev}.

A map $f: X \rightarrow Y$ between two quandles (resp. racks) $X,Y$
is called a quandle (resp. rack) homomorphism if 
$f(a*b)=f(a)*f(b)$ for any $a, b \in X$. 
A (quandle or rack) homomorphism is a (quandle or rack) isomorphism
if it is bijective. An isomorphism between  the same quandle (or rack) 
is an automorphism.

}\end{sect}

\begin{sect}{\bf Examples.\/}\label{quanxam} 
{\rm 
Any set $X$ with the operation $x*y=x$ for any $x,y \in X$ is 
a quandle called the {\it trivial} quandle. 
The trivial quandle of $n$ elements is denoted by $T_n$.

Any group $G$ 
is a quandle by conjugation as operation:
$a*b=b^{-1}ab$ for $a,b \in G$.  
Any subset of $G$  
 that is closed under conjugation is also a quandle.
For example, the set, $QS(5)$, of non-identity elements of the 
symmetric 
group on
three 
letters is a quandle. 
Similarly the set 
$$QS(6)= \{ (1234)=a, (1243)=b, (1324)=c, (1342)=B, (1423)=C,(1432)=A \}$$
of 4-cycles in the 
symmetric 
group on four letters is a quandle.

Let $n$ be a positive integer. 
For elements  $i, j \in \{ 0, 1, \ldots , n-1 \}$, define
$i\ast j= 2j-i$ where the sum on the right is reduced mod $n$. 
Then $\ast$ defines a quandle 
structure  called the {\it dihedral quandle},
 $R_n$.
This set can be identified with  the 
set of reflections of a regular $n$-gon
 with conjugation
as the quandle operation. 
We also represent the elements of $R_3$ by  $\alpha, \beta,$ and $\gamma$, 
 where the quandle 
multiplication is given by 
$x*y = z$ where $z\ne x, y$ when $x \ne y$ and $x*x=x$, for 
$x,y,z \in \{ \alpha, \beta, \gamma \}$. 

As an exercise the reader may check that 
there is a quandle homomorphism
$p:QS(6) \rightarrow R_3$ given by $p(a)=p(A)=\alpha$, $p(b)=p(B)=\beta$,
and $p(c)=p(C)=\gamma$.

Any $\Lambda={\bf Z}[T, T^{-1}]$-module $M$ is a quandle with 
$a*b=Ta+(1-T)b$, $a,b \in M$, called an {\it  Alexander  quandle}. 
Furthermore for a positive integer 
$n$, a {\it mod-$n$ Alexander  quandle}
${\bf Z}_n[T, T^{-1}]/(h(T))$
is a quandle 
for 
a Laurent polynomial $h(T)$.
The mod-$n$ Alexander quandle is finite 
if the coefficients of the  
highest and lowest degree terms 
of $h$  
 are $\pm 1$.

See \cite{Brieskorn},
 \cite{FR}, \cite{Joyce}, or \cite{Matveev} 
 for further examples.

}\end{sect}

\begin{sect}{\bf Remark.\/}
{\rm 
Let $X$ denote a quandle.
{}From Axiom~II, each element $b \in X$ defines a bijection
$S(b) : X \to X$ with $aS(b) = a \ast b$. The bijection is an automorphism
by Axiom~III.
For a word $w = b_1^{\epsilon_1} \dots b_n^{\epsilon_n}$ 
where
$b_1, \dots, b_n
\in X;
\epsilon_1, \dots, \epsilon_n \in \{\pm 1\}$,
we define
$a \ast w = aS(w)$ by
$aS(b_1)^{\epsilon_1}\dots S(b_n)^{\epsilon_n}$.
An automorphism of $X$ is called an {\it inner-automorphism\/} of $X$ if it
is $S(w)$ for a
word $w$. (The notation $S(b)$ follows Joyce's paper \cite{Joyce} and $a
\ast w$ ($= a^w$)
follows Fenn-Rourke \cite{FR}.)
}\end{sect}

\vspace{5mm} 

 Let $C_n^{\rm R}(X)$ be the free 
abelian group generated by
$n$-tuples $(x_1, \dots, x_n)$ of elements of a quandle $X$. Define a
homomorphism
$\partial_{n}: C_{n}^{\rm R}(X) \to C_{n-1}^{\rm R}(X)$ by \begin{eqnarray}
\lefteqn{
\partial_{n}(x_1, x_2, \dots, x_n) } \nonumber \\ && =
\sum_{i=2}^{n} (-1)^{i}\left[ (x_1, x_2, \dots, x_{i-1}, x_{i+1},\dots, x_n) \right.
\nonumber \\
&&
- \left. (x_1 \ast x_i, x_2 \ast x_i, \dots, x_{i-1}\ast x_i, x_{i+1}, \dots, x_n) \right]
\end{eqnarray}
for $n \geq 2$ 
and $\partial_n=0$ for 
$n \leq 1$. 
 Then
$C_\ast^{\rm R}(X)
= \{C_n^{\rm R}(X), \partial_n \}$ is a chain complex.

Let $C_n^{\rm D}(X)$ be the subset of $C_n^{\rm R}(X)$ generated
by $n$-tuples $(x_1, \dots, x_n)$
with $x_{i}=x_{i+1}$ for some $i \in \{1, \dots,n-1\}$ if $n \geq 2$;
otherwise let $C_n^{\rm D}(X)=0$. If $X$ is a quandle, then
$\partial_n(C_n^{\rm D}(X)) \subset C_{n-1}^{\rm D}(X)$ and
$C_\ast^{\rm D}(X) = \{ C_n^{\rm D}(X), \partial_n \}$ is a sub-complex of
$C_\ast^{\rm
R}(X)$. Put $C_n^{\rm Q}(X) = C_n^{\rm R}(X)/ C_n^{\rm D}(X)$ and 
$C_\ast^{\rm Q}(X) = \{ C_n^{\rm Q}(X), \partial'_n \}$,
where $\partial'_n$ is the induced homomorphism.
Henceforth, all boundary maps will be denoted by $\partial_n$.

For an abelian group $G$, define the chain and cochain complexes
\begin{eqnarray}
C_\ast^{\rm W}(X;G) = C_\ast^{\rm W}(X) \otimes G, \quad && \partial =
\partial \otimes {\rm id}; \\ C^\ast_{\rm W}(X;G) = {\rm Hom}(C_\ast^{\rm
W}(X), G), \quad
&& \delta= {\rm Hom}(\partial, {\rm id})
\end{eqnarray}
in the usual way, where ${\rm W}$ 
 $={\rm D}$, ${\rm R}$, ${\rm Q}$.

\begin{sect}{\bf Definition.\/} {\rm
The $n$\/th {\it rack homology group\/} and the $n$\/th {\it rack
cohomology group\/} \cite{FRS1}
of a rack/quandle $X$ with coefficient group $G$ are \begin{eqnarray}
H_n^{\rm R}(X;G) = H_{n}(C_\ast^{\rm R}(X;G)), \quad
H^n_{\rm R}(X;G) = H^{n}(C^\ast_{\rm R}(X;G)). \end{eqnarray}
The $n$\/th {\it degeneration homology group\/} and the $n$\/th
 {\it degeneration cohomology group\/}
 of a quandle $X$ with coefficient group $G$ are
\begin{eqnarray}
H_n^{\rm D}(X;G) = H_{n}(C_\ast^{\rm D}(X;G)), \quad
H^n_{\rm D}(X;G) = H^{n}(C^\ast_{\rm D}(X;G)). \end{eqnarray}
The $n$\/th {\it quandle homology group\/}  and the $n$\/th
{\it quandle cohomology group\/ } \cite{CJKLS} of a quandle $X$ with coefficient group $G$ are
\begin{eqnarray}
H_n^{\rm Q}(X;A) = H_{n}(C_\ast^{\rm Q}(X;G)), \quad
H^n_{\rm Q}(X;A) = H^{n}(C^\ast_{\rm Q}(X;G)). \end{eqnarray}

The homology group of a rack in the sense 
of \cite{FRS1} is $H_n^{\rm R}(X;G)$
and the cohomology of a quandle used in \cite{CJKLS} is $H^n_{\rm Q}(X;A)$.
Refer to
\cite{FRS1}, \cite{FRS2}, \cite{Flower}, \cite{Greene} for some calculations and
applications of the rack homology groups, and to \cite{CJKLS}, \cite{CJKS1}
for those of quandle cohomology groups.

\begin{sloppypar}
The cycle and boundary groups 
(resp. cocycle and coboundary groups)
are denoted by $Z_n^{\rm W}(X;G)$ and $B_n^{\rm W}(X;G)$
(resp.  $Z^n_{\rm W}(X;G)$ and $B^n_{\rm W}(X;G)$), 
 so that
$$H_n^{\rm W}(X;G) = Z_n^{\rm W}(X;G)/ B_n^{\rm W}(X;G),
\; H^n_{\rm W}(X;G) = Z^n_{\rm W}(X;G)/ B^n_{\rm W}(X;G)$$
where ${\rm W}$ is one of ${\rm D}$, ${\rm R}$, ${\rm Q}$.
We will omit the coefficient group $G$ if $G = {\bf Z}$ as usual.
\end{sloppypar}

Here we are almost exclusively interested in quandle homology or cohomology.

}
\end{sect}

\begin{sect}{\bf Example.\/}
{\rm
For $QS(6)$ (defined in \ref{quanxam}), 
we have computed using {\sc Mathematica} that 
$H^{\rm Q}_3(QS(6);{\bf Z}) = {\bf Z}_{24}.$
Similarly we have
$H^{\rm Q}_3(R_3;{\bf Z}) = {\bf Z}_{3}.$
We will illustrate in Example~\ref{qs62r3} that the 
homomorphism $p:QS(6) \rightarrow R_3$ induces a surjection
${\bf Z}_{24} \rightarrow {\bf Z}_3$.
Many other calculations are 
found in \cite{CJKLS,CJKS1,CJKS2,Greene}.

}\end{sect}

\section{Generalized 
 Knot Diagrams} \label{diagsec}

In \cite{FRS2} a general framework for colored diagrams was sketched. 
Such diagrams are useful in representing rack 
homology classes and homotopy classes of maps into the classifying space
of a rack. Here we review, amplify, and generalize
some of their constructions to the case of quandle homology classes. 
Abstract  knot diagrams are defined in \cite{KK} and
the relation to virtual knots \cite{Lou} is established in \cite{KK}.
We extend their definition to include 
arcs and surfaces with boundary. 
First we generalize crossings of classical knot diagram to higher
dimensions.

\begin{figure}
\begin{center}
\mbox{
\epsfxsize=3in 
\epsfbox{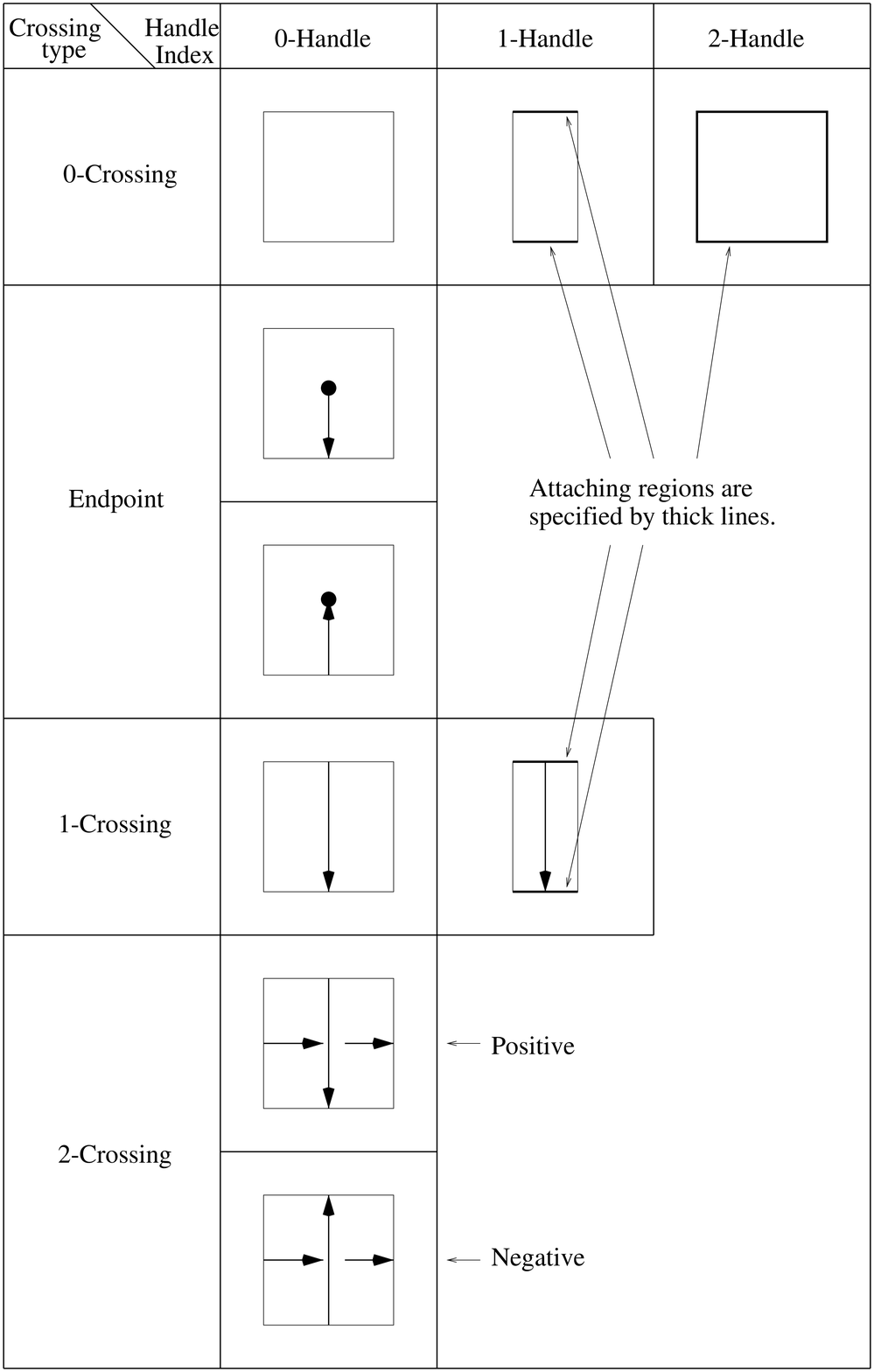}
}
\end{center}
\caption{A table of handles for abstract 1-knot diagrams }
\label{handlesindim2}
\end{figure}

\begin{figure}
\begin{center}
\mbox{
\epsfxsize=4in
\epsfbox{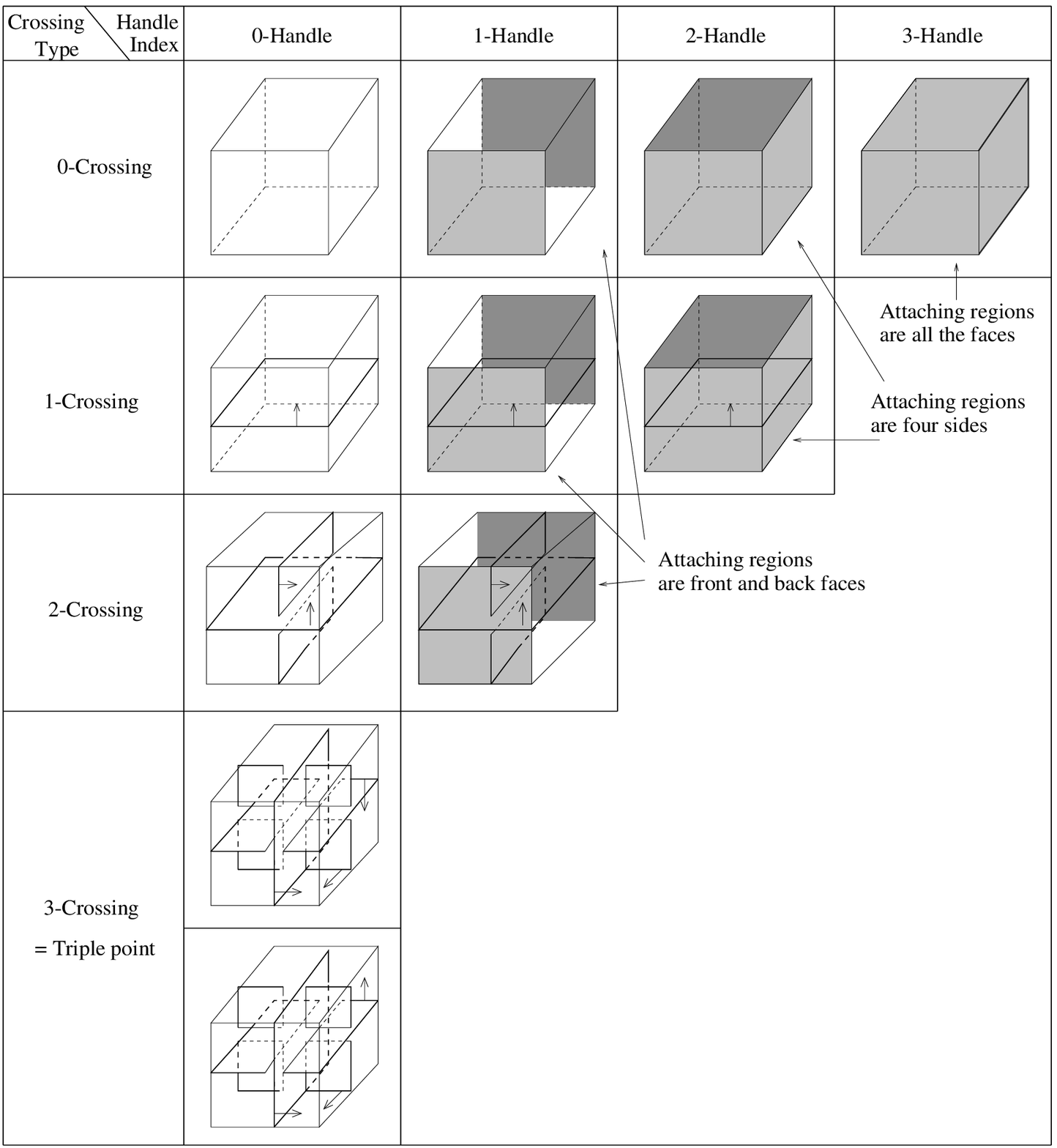}
}
\end{center}
\caption{A table of handles for abstract 2-knot diagrams, crossings }
\label{handlesindim3}
\end{figure}

In the top 
row  
and the 
bottom three rows of 
Fig.~\ref{handlesindim2}, 
 $k$-crossing $1$-diagrams
are indicated for 
$k=0,1,2$. In Fig.~\ref{handlesindim3}
$k$-crossing $2$-diagrams are indicated for $k=0,1,2,$ and $3$. 
We define such diagrams  
in general.

\begin{sect} {\bf Definition.\/} 
\label{kndiag}
{\rm 
Fix a non-negative integer $n$.
Let $E$ be the $(n+1)$-disk $[-1,1]^{n+1}
=\{ (x_1, \dots, x_{n+1})
| -1 \leq x_i \leq 1, i=1,  \cdots, n+1 \}$, 
 and for $j\in \{1,\dots, n+1\}$,
let $E_j$
be the $n$-disk in $E$ determined by $x_j=0$.
A {\it $0$-crossing $n$-diagram} is the disk $E$.
A {\it $1$-crossing $n$-diagram} (or simply 1-crossing diagram)
$\Sigma_1$ 
is the 
pair of $E$ and the 
$n$-disk $E_{1}$.
Let $k$ be an integer in the range $\{1, \ldots , n+1\}$.
A {\it $k$-crossing $n$-diagram} (or simply $k$-crossing diagram)
$\Sigma_k$ 
is the pair of $E$ and the
union of the disks 
$E_{1}\cup \cdots  \cup E_{k}$ with crossing information
that is indicated by removing  
from 
$E_{j}$ (for $j=2, \ldots, k$) an open $\varepsilon$-neighborhood 
of $(E_{1} \cup \cdots \cup E_{j-1}) \cap E_{j}$.
For example, when $n=1$ a $2$-crossing $1$-diagram is the standard 
depiction of a classical knot crossing. When $n=2$, a $2$-crossing $2$-diagram
is the broken surface depiction \cite{CS:book} 
 of pair of planes in ${\bf R}^4$
that project to intersect along a double curve
as depicted in Fig.~\ref{handlesindim3} third entry from top. 
 A $3$-crossing $2$-diagram is a 
broken diagram of a triple point 
as depicted in  Fig.~\ref{handlesindim3} bottom two entries.

The removal of the $\varepsilon$-neighborhood is a convention for
indicating relative height in a projection from $(n+2)$-dimensional space.
Thus in $E\times [0,1]$, we re-embed  $E_k$  as 
$\tilde{E_k} = E_k \times \{ \frac{n+2-k}{n+2} \}$. 
Then {\it the lift},
$\tilde{E_k}$,
 projects to $E_k$ under the projection
$E \times [0,1] \rightarrow E$. The relative position in the
$[0,1]$ factor then is determined by the number of regions on the
 sheet
$E_k$.

Such diagrams, $(E,\cup E_j)$, are oriented.
For example, when $n=1$, orientations of arcs are indicated 
by arrows along the arc. 
The  orientation
of $E$ is indicated by a short arrow  normal to the arc $E_j$
where tangent followed by normal is the standard 
orientation of the disk.
The normal arrow may be omitted 
from the figure when the orientation of the
disk $E$ is that of the plane of the paper.
When $n=2$, we assume that the disks $E_j$ are oriented in a counterclockwise fashion and only indicate the normal arrow.

In higher dimensions, 
signs of $(n+1)$-crossing $n$-diagrams
 are determined by normal vectors as follows.
Let the normal to $E_j$ be denoted by $\nu_j$. Then
the $(n+1)$-tuple 
$(\nu_1, \ldots, \nu_{n+1})$ determines an orientation of
the disk $E$. The $(n+1)$-crossing point is {\it positive} 
if and only if
 this orientation agrees with the standard orientation. 
Thus if the normal to $E_j$ is the  
standard basis vector $e_j =(0,\ldots , 1, \ldots, 0)$
(where the non-zero entry is in the $j\/$th 
coordinate), then
the $(n+1)$-crossing point is positive.
}\end{sect}

\begin{sect}
{\bf  Definition.\/}
{\rm There are $2^k-1 = 1+ 2+ \cdots +2^{k-1}$ components in $\Sigma_k$.
We call them the  {\it $n$-regions\/} of $\Sigma_k$.
The disk $E_j$ is called the {\it level $j$ sheet\/} of $\Sigma_k$.
The {\it $k$-crossing point set} (or simply $k$-crossing) in $E$
is the intersection 
$E_{1} \cap E_2 \cap \cdots \cap E_{k}$,
which is
a disk of dimension $n+1-k$ (or of codimension $k$).
The complement of $E_{1}\cup \cdots \cup E_{k}$ 
in $E$ has $2^k$ components called 
{\it $(n+1)$-regions}.
}\end{sect}

\begin{sect}{\bf  Remark.\/} \label{subdiag} 
{\rm 
Observe that in the level $j$ sheet at a 
$k$-crossing we have a  
$(j-1)$-crossing diagram; this diagram is the 
orthogonal
projection of the previous levels onto the $j\/$th level. 
}\end{sect}

Next we introduce 
endpoint,
branch point, 
and hem diagrams.
These exceptional diagrams will be useful to us in describing
degenerate chains in quandle homology.

\begin{sect}{\bf Defintion.\/} {\rm  
An {\it endpoint diagram} is a 
$2$-disk with an 
arc embedded; one end of the arc is in the interior of the disk, 
the other end is on the boundary, and the arc intersects
 the boundary transversely. In our illustrations, we indicate
 the interior point as a fat vertex. 
 See Fig.~\ref{handlesindim2} the second and third pictures from top. 
}\end{sect}

\begin{figure}
\begin{center}
\mbox{
\epsfxsize=4in
\epsfbox{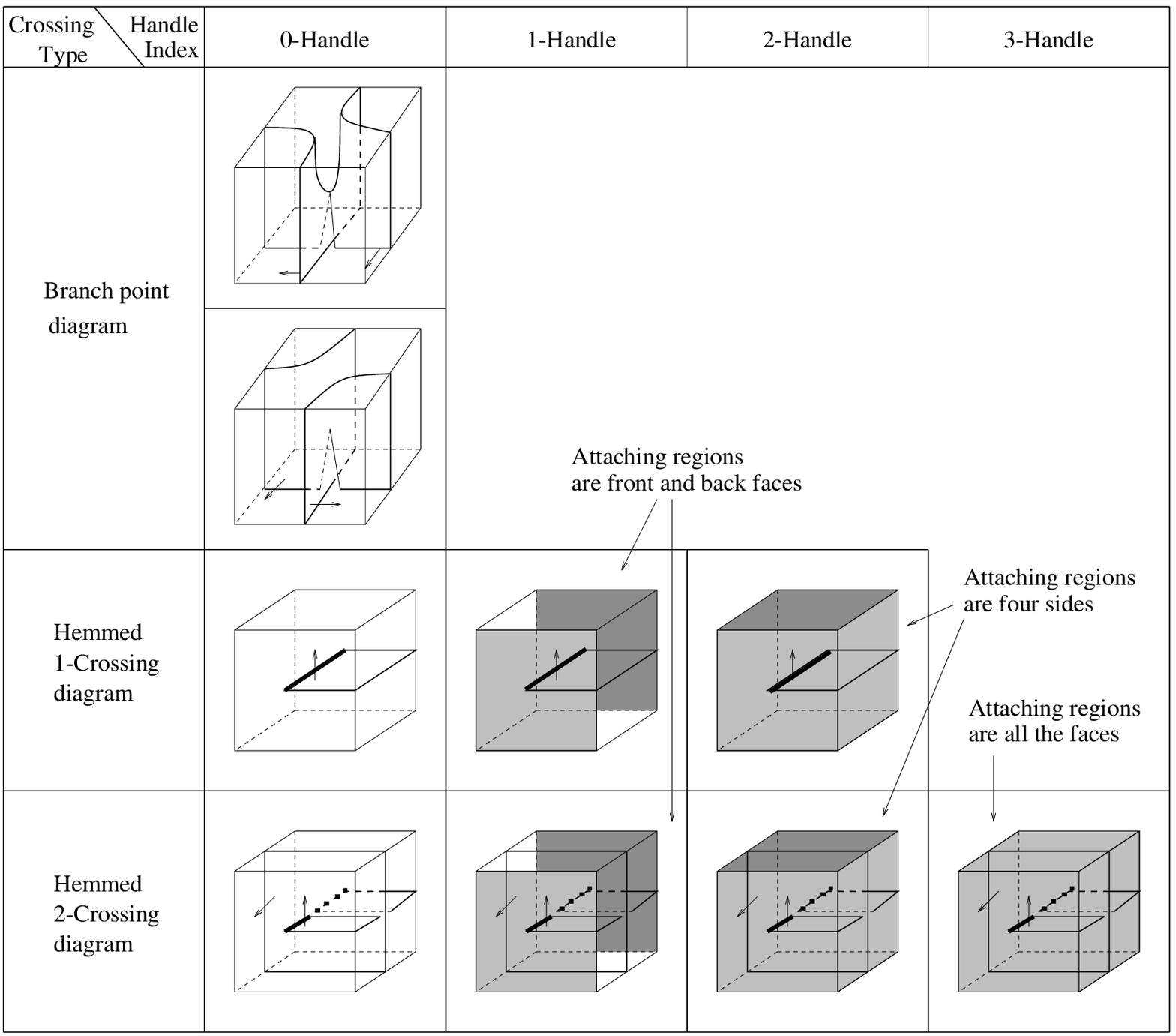}
}
\end{center}
\caption{A table of handles for abstract 2-knot diagrams, branch points and 
end curves }
\label{handlesindim3II}
\end{figure}

\begin{sect}{\bf Definition.\/}
{\rm A {\it  branch point diagram} consists of a 
neighborhood of 
a branch point (also called Whitney umbrella, a generic 
singularity 
of surface maps) in 
a $3$-ball
in which over and under crossing information is indicated
by removing
a thin open triangle 
along the double point arc
as depicted in
Figure~\ref{handlesindim3II}, the top and second top pictures.
The vertex of the triangle is at the branch point
 which is in the interior of the 3-ball.
The boundary of the branch point neighborhood on the 
boundary of the $3$-ball is the  diagram of the unknot 
with one crossing.
 }\end{sect}

\begin{sect} {\bf Definition.\/}
{\rm A 
{\it hemmed $1$-crossing diagram} 
is a $2$-dimensional half-disk embedded in
a $3$-ball as depicted in Fig.~\ref{handlesindim3II} 
the second bottom pictures. 
A part of the boundary of the disk is in the interior of the $3$-ball, 
called the {\it hem}. 

} \end{sect}

\begin{sect} {\bf Definition.\/}
{\rm A {\it 
hemmed $2$-crossing diagram}
consists of 
two disks immersed (with crossing information) in the 3-ball with the boundary of one of these
(the under-sheet) contained in the interior of the $3$-ball.
Crossing information is depicted in the bottom
of Fig.~\ref{handlesindim3II} by removing a thin rectangular 
neighborhood of the double point arc in the under-sheet ---
the horizontal sheet in the illustration. 
The under-sheet of a hem $2$-crossing diagram 
is the sheet that has an interior 
hem boundary (depicted by a thick line segment).
It is important to note that 
 at a hemmed $2$-crossing diagram
it is the hemmed surface that is 
the under-sheet.

 } \end{sect}

An endpoint diagram, a branch point diagram, or
a   hemmed $1$- or 
 $2$-crossing diagram 
is also 
 called an {\it exceptional diagram}.

In Figs.\ref{handlesindim2}, \ref{handlesindim3}, and \ref{handlesindim3II},
the collection of endpoint diagrams, 
branch point diagrams, 
crossing diagrams 
in low dimensions,
and the handles that they 
correspond 
to are indicated.

\begin{figure}
\begin{center}
\mbox{
\epsfxsize=1.5in
\epsfbox{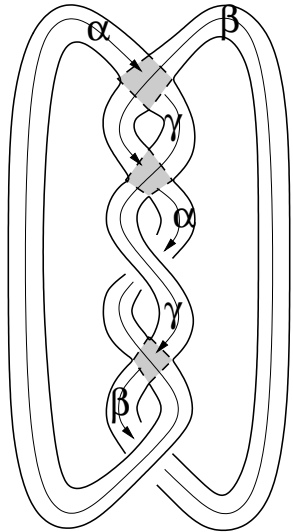} 
}
\end{center}
\caption{A $1$-knot diagram }
\label{ssvsv1}
\end{figure}

\begin{figure}
\begin{center}
\mbox{
\epsfxsize=1.7in
\epsfbox{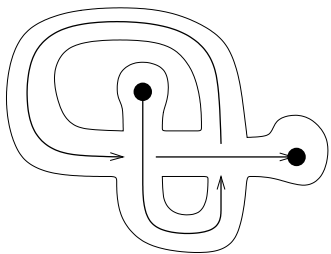} 
}
\end{center}
\caption{A $1$-knot (arc) diagram }
\label{arcdiag}
\end{figure}

\begin{figure}
\begin{center}
\mbox{
\epsfxsize=5in  
\epsfbox{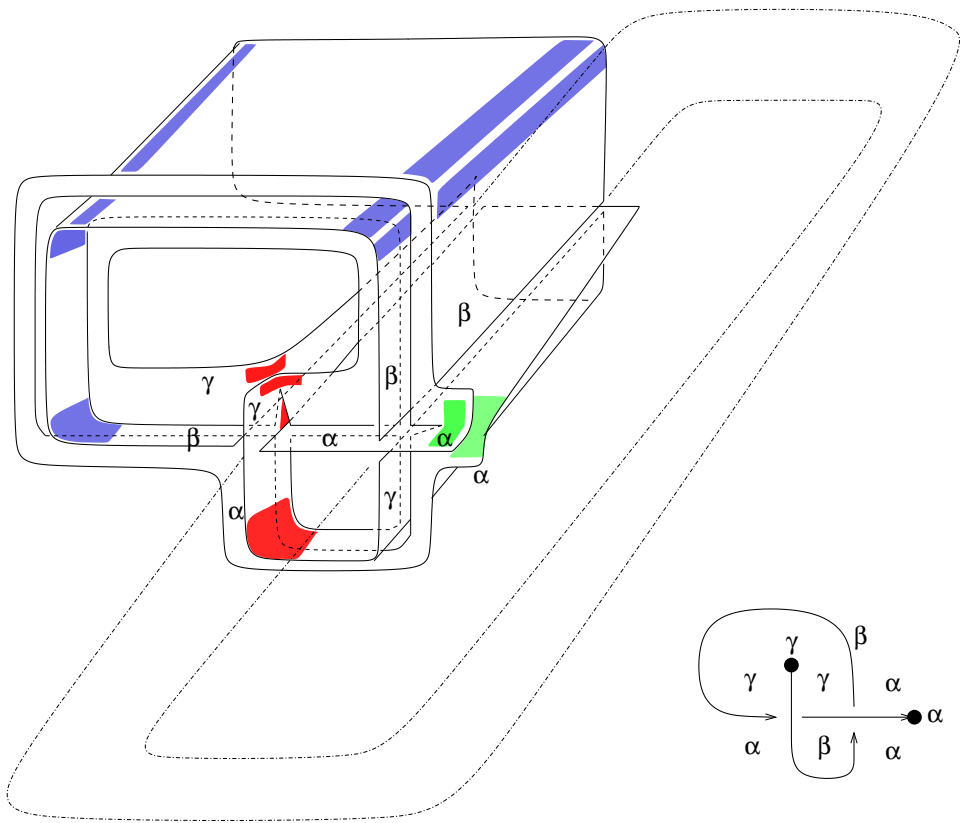} 
}
\end{center}
\caption{An abstract $2$-knot diagram }
\label{genh3r3susp}
\end{figure}

\begin{sect} {\bf Definition.\/} \label{absdef} 
{\rm 
A
{\it generalized 
abstract 
 $0$-knot diagram} 
consists of a collection of 
vertices (points) 
embedded in a closed $1$-manifold.
Let $n=1$ or $n=2$. 
A
generalized 
abstract $n$-knot diagram is 
the image of  
a 
 generic map $f:M\rightarrow N$
with crossing information indicated
where:  
\begin{enumerate}

\item
$f(M)\cap \partial N= \emptyset$;

\item 
Local crossing information indicated as follows:
Each point of $N$ has a neighborhood $B$, such that
$B\cap f(M)$ is a $k$-crossing $n$-diagram ($k=0,\ldots, n+1$),
or an  exceptional diagram.
\end{enumerate}
The image 
$f(M)$ is called the
{\it universe of the diagram}. If the manifold $M$ is closed, then the diagram is said to be {\it closed.} Thus when $n=2$, a closed diagram 
may have branch points but does not have hems.

Abstract diagrams will often be denoted by 
$K=[f:M\rightarrow N]$. 
This is a slight abuse of notation since we are considering the image 
$f(M)\subset N$ with crossing information.

Abstract diagrams 
are depicted in Figs.~\ref{ssvsv1}, 
\ref{arcdiag},
\ref{genh3r3susp}, 
and \ref{example0}.

}\end{sect}

\begin{sect}{\bf Remark.\/}
{\rm
In \cite{NK} abstract link diagrams were introduced in which 
the image $f(M)$ is a deformation 
retract of the ambient space $N$ (see also
\cite{KK}).  Here we do not need 
to use this strong condition. Moreover,
a generalized
abstract diagram (of a closed manifold $M$) in the current sense 
gives rise to an  abstract diagram 
in the sense of \cite{NK} 
by taking a regular neighborhood of the image. 
Figure~\ref{arcdiag} is 
a strict generalization of 
abstract diagrams defined in \cite{NK}, 
and Fig.~\ref{ssvsv1}  is an abstract diagram in any sense.
Throughout the sequel, 
generalized abstract diagrams are called abstract diagrams
for simplicity.
} \end{sect}

Abstract 
diagrams can be constructed via handle theoretic techniques.
In Figs.~\ref{handlesindim2}, \ref{handlesindim3},  and 
\ref{handlesindim3II} 
crossing diagrams 
and 
exceptional diagrams are classified as 
handles,
 and 
attaching 
regions are indicated.

\begin{sect} {\bf Remark.\/} {\rm
Let $K$ be an abstract $1$-knot diagram in $N$. 
Then $K$ determines an embedded arc $\hat{K}$ in $N \times [0,1]$
by lifting  each 
over-arc to a higher 
level (in the 
$[0,1]$ factor)
than the under-arc at each double point of $K$,
as is  typically done in classical knot theory.
In the present convention endpoints are embedded in
$N\times \{0\}$.
} \end{sect}

\begin{sect}{\bf Example.\/}
{\rm
An abstract diagram of a closed curve is 
indicated in Fig.~\ref{ssvsv1}.  
The $1$-manifold  is $S^1$, and the manifold $N$ is a 
thrice 
punctured 
torus. There are three $2$-crossing $1$-diagrams as its subdiagrams, 
the top two  
are 
positive  and 
the bottom one is negative.
The $0$-handles of $N$ 
are these 
$2$-crossings 
 depicted as lightly 
shaded squares.
The unshaded bands are $1$-handles that 
are 
$1$-crossing 
$1$-diagrams.

}\end{sect}

\begin{sect}{\bf Example.\/}
{\rm
An abstract diagram of an arc in 
a planar surface 
is indicated 
 in Fig.~\ref{arcdiag}. 
The diagram has two positive crossing points and two endpoints.

}\end{sect}

\begin{sect} {\bf Remark.\/} {\rm 
An abstract 
closed  
 $2$-knot diagram is constructed from copies of $k$-crossing 
$2$-diagrams ($k=1,2,3$) and branch point diagrams as follows.
The crossing diagrams are identified along their boundaries so that
the arcs of double points either form closed loops or end at 
$3$-crossing diagrams or branch point diagrams. Thus in the boundary of 
$N$, there is, at most,
 a collection of simple closed curves.  
These then are capped-off by adding $2$-handles of the form
$1$-crossing $2$-diagrams.
Moreover, when crossing 
diagrams agree at points of their boundaries the 
broken and unbroken sheets 
match broken and unbroken sheets, respectively
(see Fig.~\ref{gluecond}, 
Fig.~\ref{handlesindim3} bottom two pictures, 
and Fig.~\ref{duncehat}). 
It is also required that orientations match when crossing diagrams are 
glued.

\begin{figure}
\begin{center}
\mbox{
\epsfxsize=2.5in
\epsfbox{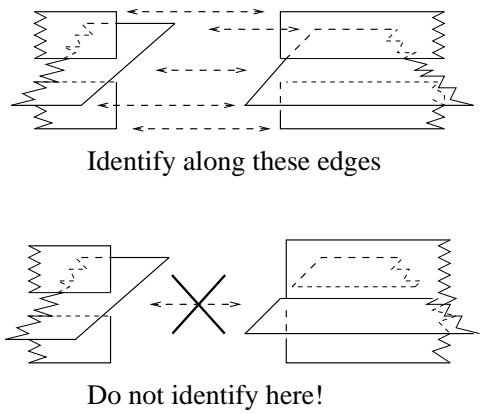} 
}
\end{center}
\caption{Broken sheets are glued together consistently} 
\label{gluecond}
\end{figure}

} \end{sect}

\begin{sect} {\bf Remark.\/} {\rm 
An abstract 
$2$-knot diagram represents 
a surface $M$ embedded 
in $N \times 
[0,1]$ by filling in the surface along the broken sheets in the $[0,1]$ 
direction (see \cite{CS:book} for details on filling in broken surface 
diagrams). In this way,
 the surface represented by the diagram is 
embedded in $N\times [0,1]$.

} \end{sect}

\begin{sect}{\bf Example.\/}
{\rm In Fig.~\ref{genh3r3susp}, 
a part 
of a 
closed 
abstract $2$-knot 
diagram is indicated (the 
 relation to the arc diagram 
in the 
lower right
will be explained 
subsequently). The diagram has two 
$3$-crossing points (triple points) 
and two branch points. The arcs of 
double points that end in the 
front of the diagram  are to be glued to 
the back of the diagram via  
two $1$-handles each of which 
is a 
$2$-crossing diagram.  
 In this way, 
the $3$-manifold that contains the 
diagram is a handle-body, and the 
boundary contains 5 simple closed curves. 
These 
are 
attaching regions 
for $2$-handles that 
are
$1$-crossing diagrams.

}\end{sect}

\begin{sect} {\bf Remark.\/}
{\rm Abstract diagrams can be defined similarly in higher dimensions
using $k$-crossing $n$-diagrams, as in \cite{FRS2}.
However, the branch point set is 
more  
subtle in higher dimensions. 
} \end{sect}

\section{Coloring Abstract 
Diagrams} \label{colorsec}

Abstract diagrams, when colored by 
quandles, represent homology classes. 
Colored abstract $1$- and $2$-knot diagrams 
are of particular interest.

\begin{figure}
\begin{center}
\mbox{
\epsfxsize=4in
\epsfbox{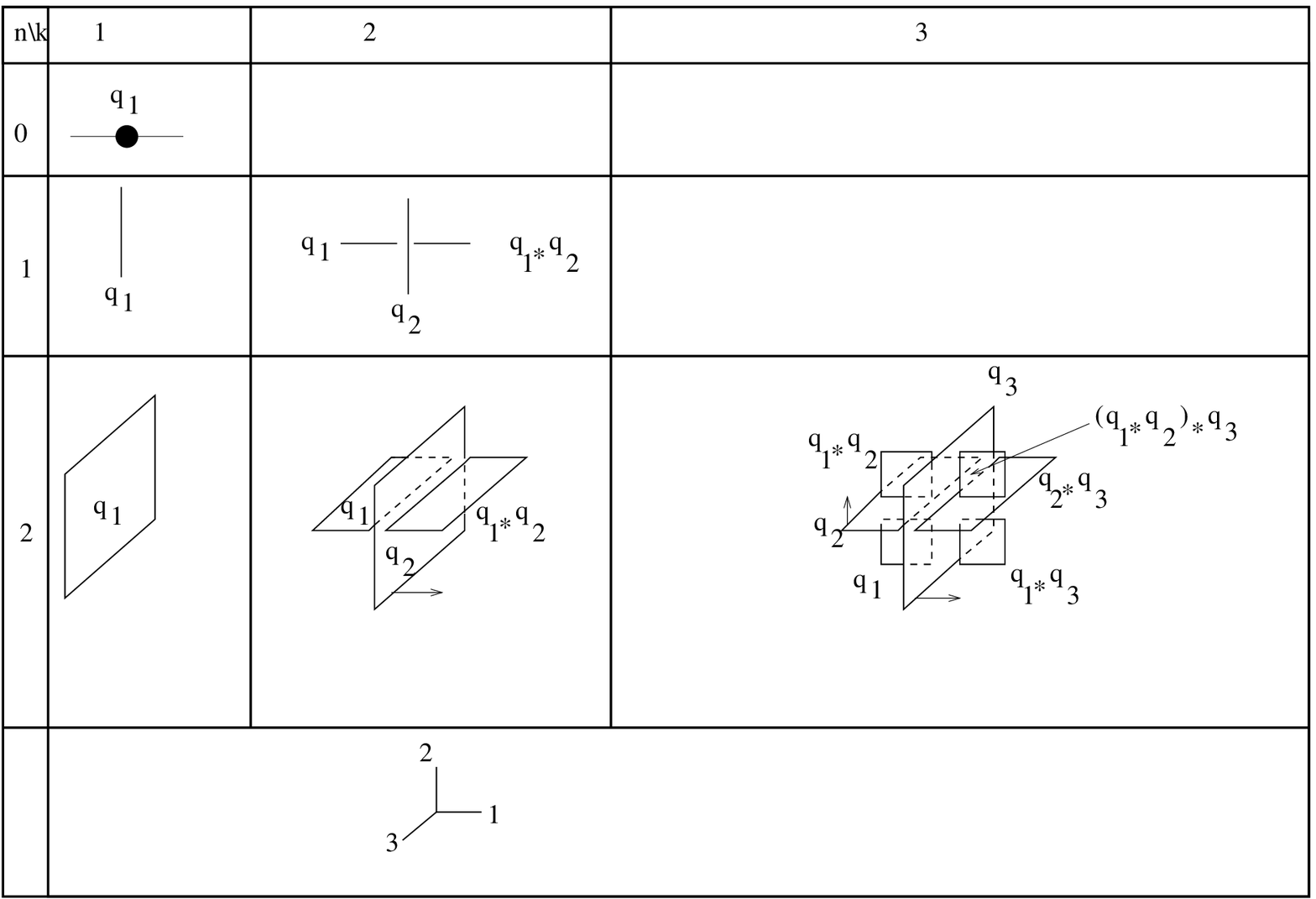}
}
\end{center}
\caption{Colors of $k$-crossing $n$-diagrams }
\label{kcrossingndiag1}
\end{figure}

\begin{sect}{\bf Definition.\/}
{\rm  
For $n=0,1,2$ and $k\le n+1$, 
{\it a coloring of a $k$-crossing $n$-diagram} 
is defined as an assignment of quandle 
elements to the points, arcs or regions of 
the diagram such that the following 
{\it color condition at the crossing} is satisfied.
The condition is depicted in 
Fig.~\ref{kcrossingndiag1}. In this figure,
an under-arc (or lower sheet of a 
$2$-dimensional region) is colored 
with a quandle element, say  $q_1$. The arc 
(or region) 
that is so colored is the arc away 
from which the normal to the over-arc
(or region) points. The over-arc is 
colored with $q_2$. And the remaining 
under arc (or region) is colored with 
$q_1 * q_2$. In case $n=0$ only the 
point is colored with $q_1$. In case 
the crossing point is a triple point, 
then the middle sheet has 
two colors ($q_2$ and $q_2 * q_3$)
as described above, and the 
lowest sheet 
has four colors ($q_1$, $q_1*q_3$, 
$q_1*q_2$, and $(q_1*q_2)*q_3$).

A 
branch point diagram is 
colored by assigning a single quandle 
element to the 2-dimensional region of the diagram.

Let $n=1,2$.
{\it A  (quandle) coloring of 
a closed 
abstract $n$-knot diagram} 
is an assignment of quandle elements 
to the connected regions of the diagram such that 
the restriction to each  
branch point (when $n=2$)
or
$k$-crossing diagram is a coloring.
More precisely, let ${\cal R}$ be the set of arcs (for $n=1$) 
or regions (for $n=2$) of 
an
abstract $n$-knot 
diagram ($n=1,2$). 
A coloring is a map ${\cal C} : {\cal R} \rightarrow X$,
where $X$ is a 
quandle, such that 
$\{ {\cal C}(r) | r \in {\cal R} \}$ 
satisfies 
the above mentioned conditions 
at each $k$-crossing $n$-diagram in 
the given abstract $n$-knot diagram.

}\end{sect}

\begin{figure}
\begin{center}
\mbox{
\epsfxsize=5in
\epsfbox{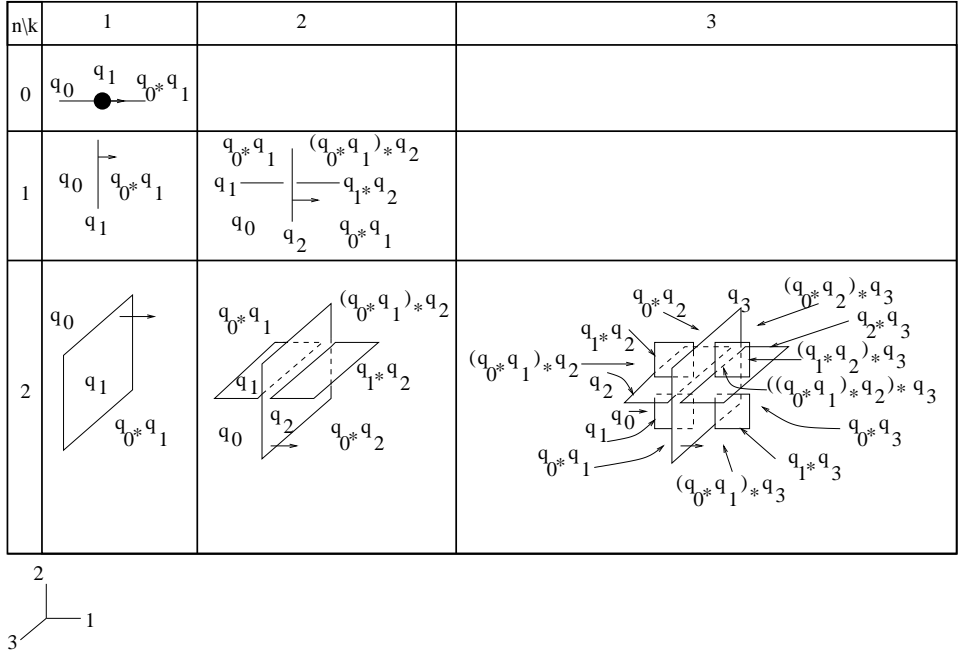}  
}
\end{center}
\caption{Shadow colors of $k$-crossing $n$-diagrams }
\label{kcrossingshad} 
\end{figure}

\begin{sect} {\bf Examples.\/} 
{\rm  Figure~\ref{ssvsv1} indicates a coloring of 
an abstract $1$-knot diagram by $R_3$.
Figure~\ref{genh3r3susp} indicates a coloring of 
an abstract $2$-knot diagram by $R_3$.
} \end{sect}

In Fig.~\ref{kcrossingshad}     
shadow colorings 
of $k$-crossing 
$n$-diagrams
for all  $n= 0,
1,2,$ and $1\le k\le n+1$ 
are depicted. We 
define
these case by case.

\begin{sect}{\bf Definition.\/}
{\rm 
For any $n$, 
{\it a shadow coloring of a 
$0$-crossing diagram} is an assignment of a quandle 
element to the $(n+1)$-dimensional region of the diagram.
{\it A  shadow coloring of a $1$-crossing 
diagram} is an assignment of three 
quandle elements to the constituents of 
the diagram, as follows. 
The $(n+1)$-dimensional region away from  which 
the normal of the  $1$-crossing recieves a quandle element $q_0$, 
the $n$-dimensional region ($1$-crossing) receives the color $q_1$,
and the $(n+1)$-dimensional region into which the normal
points receives 
color $q_0*q_1$.

In general, {\it a shadow coloring of a
$k$-crossing $n$-diagram} 
is an assigment of colors 
to the $n$ and $(n+1)$-dimensional regions
that satisfies the following conditions. 
Along the $n$-dimensional regions of the diagram, 
it is a coloring of the 
$k$-crossing diagram.  
Each point in an $(n+1)$-dimensional 
region of the diagram is assigned a quandle 
element, and if two such regions are separated 
by an $n$-dimensional region, then  any point in the 
$n$-dimensional region has a neighborhood 
homeomorphic to a $1$-crossing $n$-diagram. In this case, 
the coloring on the $(n+1)$-dimensional regions 
satisfies the condition that the restriction 
to such a neighborhood is a 
shadow coloring of the $1$-crossing diagram.

\begin{figure}
\begin{center}
\mbox{
\epsfxsize=6in
\epsfbox{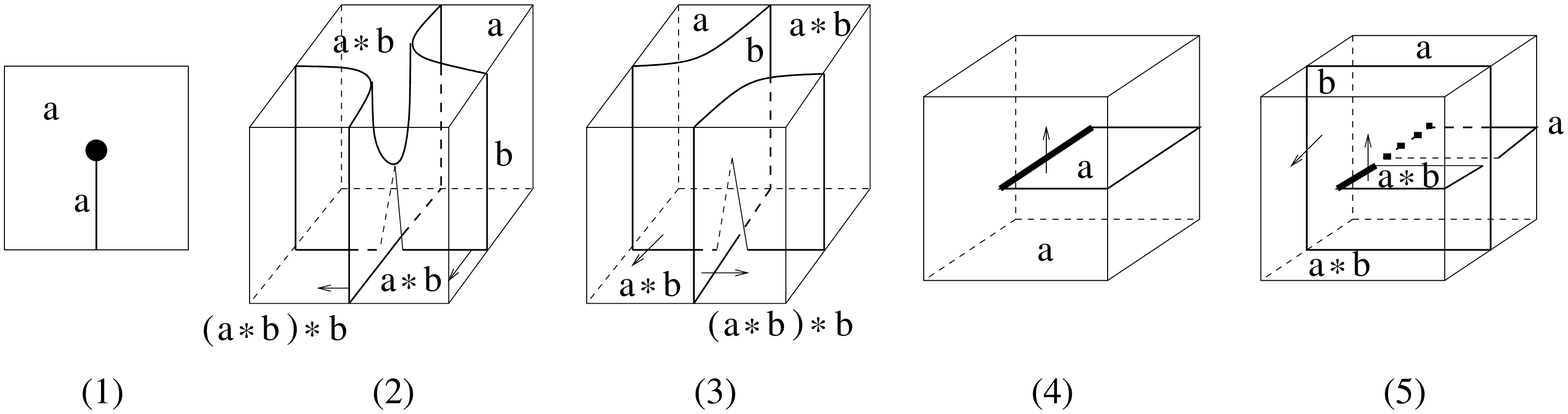}
}
\end{center}
\caption{Shdow colors of end/branch points and hems }
\label{shadows}
\end{figure}

{\it A shadow coloring of an 
endpoint diagram} is the assignment of 
the same color to the arc and the 
adjacent region. See Fig.~\ref{shadows} (1). 
A {\it shadow coloring 
of a branch point diagram} is the assignment 
of a single quandle element, say $b$, to the 
$2$-dimensional faces of the diagram. 
There are three $3$-dimensional regions in a 
neighborhood of the branch point. One such 
region is assigned the color $a$; the other 
two are assigned the color $a*b$ or $(a*b)*b$, 
and the assignment is determined by the 
condition that normal vectors point 
into regions that receive quandle products.  
See Fig.~\ref{shadows} 
(2) and (3).

A {\it shadow coloring 
of a hem $1$-crossing diagram} 
is an assignment of a single quandle element 
to the region and the disk.  See Fig.~\ref{shadows} (4).
A {\it shadow coloring 
of a hem $2$-crossing diagram} is an assignment of
an element $a$ on the disk and the region away from which the normal 
of the disk dividing the hem, the dividing disk receives $b$, 
and then, the region and the disk into which the normal 
points receive the element $a*b$. 
See Fig.~\ref{shadows} (5).

Let $n=1,2$.
{\it A  shadow coloring of an abstract $n$-knot diagram} 
is an assignment of quandle elements 
to the $n$  and $(n+1)$-dimensional  regions of the diagram such that 
the restriction to each  endpoint, branch point, hem 
 or 
$k$-crossing diagram is a 
shadow coloring.

A shadow coloring is, as before, 
 a map ${\cal C}: {\cal R} \rightarrow X$ 
where   ${\cal R}$ is the set of
arcs ($n=1$) or regions ($n=2$) 
and complementary $(n+1)$-dimensional regions.
}\end{sect}

\begin{figure}
\begin{center}
\mbox{
\epsfxsize=2in
\epsfbox{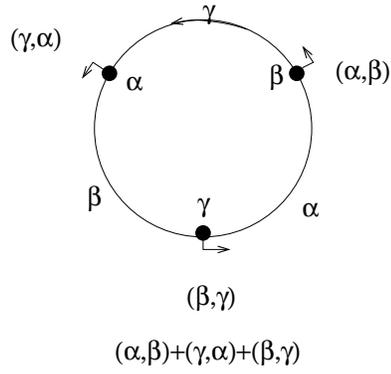}
}
\end{center}
\caption{A  shadow colored abstract $0$-knot  diagram }
\label{example0}
\end{figure}

\begin{sect}{\bf Example.\/}
{\rm Figure~\ref{example0} indicates a shadow quandle coloring 
by $R_3$ of an abstract $0$-knot diagram. 
The small flag-like arrows at the 1-crossing points indicate the 
orientation of the vertices. 
 When this arrow coincides with 
the tangent direction of the circle (as in the figure), 
the crossing is positive; otherwise the crossing is negative.
} \end{sect}

\begin{sect}{\bf Example.\/}
{\rm Figure~\ref{genh3r3} indicates a shadow quandle coloring 
by $R_3$ of the 
abstract $1$-knot diagram
that was 
given in Fig.~\ref{arcdiag}. 
The boundary of the surface is abbreviated and is not drawn in
Fig.~\ref{genh3r3} for simplicity. 
In the figure, colors on regions are indicated by letters in squares.

\begin{figure}
\begin{center}
\mbox{
\epsfxsize=1.7in
\epsfbox{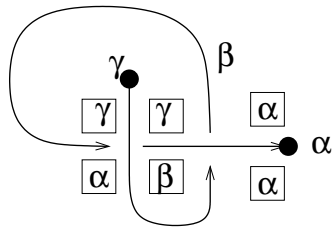}
}
\end{center}
\caption{A  shadow colored abstract $1$-knot 
diagram }
\label{genh3r3}
\end{figure}

}
\end{sect}

\begin{sect}{\bf Remark.\/}
{\rm 
Shadow colorings were  defined in \cite{FRS2}
and used in \cite{RS} to prove that the right and left 
trefoil knots are distinct.
The notions of coloring 
and shadow coloring of $k$-crossing $n$-diagrams 
extend for all $n=1,2, \ldots,$ and for all 
$k=0, \ldots, n+1$ (see 
 \cite{FRS2}). 
}\end{sect}

To see that shadow colorings 
exist, we prove the following lemmas.

\begin{sect} {\bf Lemma.\/} \label{colorlemma}
Let ${\cal C}$ be a coloring 
by a  quandle $X$ 
of an oriented $k$-crossing $n$-diagram, 
$(k \leq n)$,  
which is 
obtained from 
an oriented $(k+1)$-crossing $n$-diagram 
by ingoring the level $k+1$ sheet, $E_{k+1}$.
Let $D$ be an $n$-region  
on 
$E_{k+1}$,
and let $x \in X$. 
Then there is a unique 
coloring  
${\cal C}'$ by $X$ of the $(k+1)$-crossing diagram
which restricts to ${\cal C}$ and such that ${\cal C}'(D)=x$.
\end{sect}
{\it Proof.\/}
Let $D_j$, $j=1, \cdots, 2^{k+1}$, be the $n$-regions 
on
$E_{k+1}$.
Let $D=D_1$, and pick points $p_j \in D_j$. 
For any fixed $j$, let $\gamma$  be  a path in the level 
$k+1$ sheet 
 (the image of a continuous map 
from $[0,1]$ to the level  $k+1$ sheet) 
 from $p_1$ to $p_j$. 

Assume without loss of generality 
that the path is in general position
with the level $g$ sheets for 
all $g$, so that it meets the sheets 
in  finitely many points. Let $s_1, \cdots, s_r$ 
be the intersection points,
and let $\epsilon_i$, $i=1, \cdots, r$ be $+1$ (resp. $-1$)
if the path goes in the same (resp. the opposite) direction as
the normal of the sheet at $p_i$. 
Let ${\cal C}(p_i)=c_i$ 
(precisely speaking, the color of the disk in which 
the point $p_i$ lies). 
Then define the color ${\cal C}'(D_j)=x * w$
where $w$ is the word $ c_1^{\epsilon_1} c_2^{\epsilon_2} 
\cdots c_n^{\epsilon_r}.$
This definition is made in such a 
way that the condition of coloring
is satisfied 
at each intersection point along the path  $\gamma$, and 
in remains to be proved that the color thus defined does not 
depend on the choice of the path $\gamma$. 

Let $\gamma_i$, $i=0,1$, be such two paths. 
Since the level $k+1$ sheet 
is simply connected,
 there is a homotopy between them.
Such a homotopy is a map from a $2$-disk to the level 
$k+1$ sheet, 
whose image is denoted by $U$.
Assume without loss of generality that $U$ is in general position 
with the level $g$ sheets for all $g (<k)$. 
Then the intersection between $U$ and the level $g$ sheets for all $g (<k)$
is generically immersed $1$-manifold with boundary, i.e., arcs 
with transverse double points. When the 
$\varepsilon$-neighborhood 
 in the definition of the crossing diagram is removed from 
these arcs, then we obtain a classical knot diagram on $U$ 
with boundary of arcs liying on $\partial U$. 
Give the color by $X$ on the regions and arcs in $U$ by 
using the rule of the shadow color depicted in 
Fig.~\ref{shadows}.
Then the color is well-defined on the diagram on $U$, 
proving that the color at the end point of $\gamma_i$, $i=0,1$, 
coincide.
$\Box$

\begin{sect} {\bf Lemma.\/} \label{1connlemma}
Let 
$K=[f:M^n\rightarrow N^{n+1}]$ 
be an
abstract $n$-knot diagram,
and  let ${\cal C}$ be a 
coloring 
by a  quandle $X$.
Suppose $N$ is simply connected.
Then for any 
$(n+1)$-region $R$
and $x \in X$, 
there is a shadow color ${\cal C}'$ 
for 
$K$ 
such that  ${\cal C}'$ 
restricts  ${\cal C}$ and 
  ${\cal C}'(R)=x$. 
\end{sect} 
{\it Proof.\/}
This is proved by defining 
colors using paths as in Lemma~\ref{colorlemma},
and the proof is similar. 
The condition $\pi_1(N)=\{1\}$ guarantees
 the existence of a homotopy between paths.

In the case $n=1$ and $2$, we included end/branch points
and hem $1$- and $2$-crossings, 
in addition 
to $k$-crossing $n$-diagrams.
For these diagrams, a 
similar argument applies, using the definition 
of shadow colors for these diagrams.
$\Box$

\begin{figure}
\begin{center}
\mbox{
\epsfxsize=5in
\epsfbox{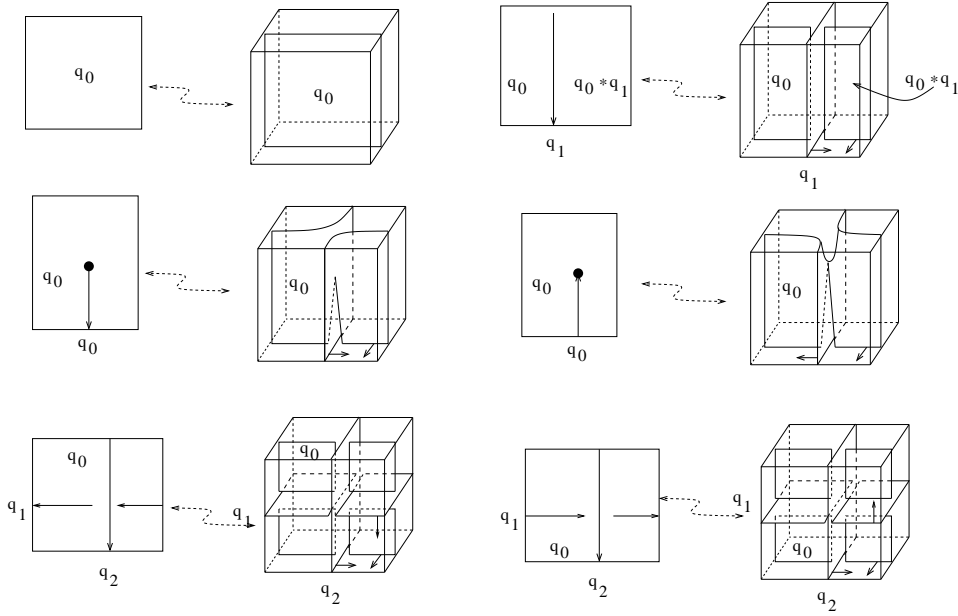}
}
\end{center}
\caption{The correspondence between shadow 
colored 1-knot diagrams and colored 2-knot diagrams }
\label{ab1eq2}
\end{figure}

\begin{figure}
\begin{center}
\mbox{
\epsfxsize=3.5in
\epsfbox{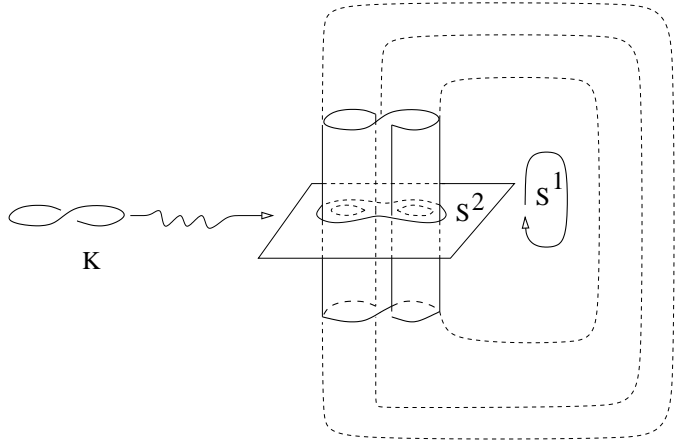}
}
\end{center}
\caption{The suspension construction }
\label{susp}
\end{figure}

\begin{figure}
\begin{center}
\mbox{
\epsfxsize=5in
\epsfbox{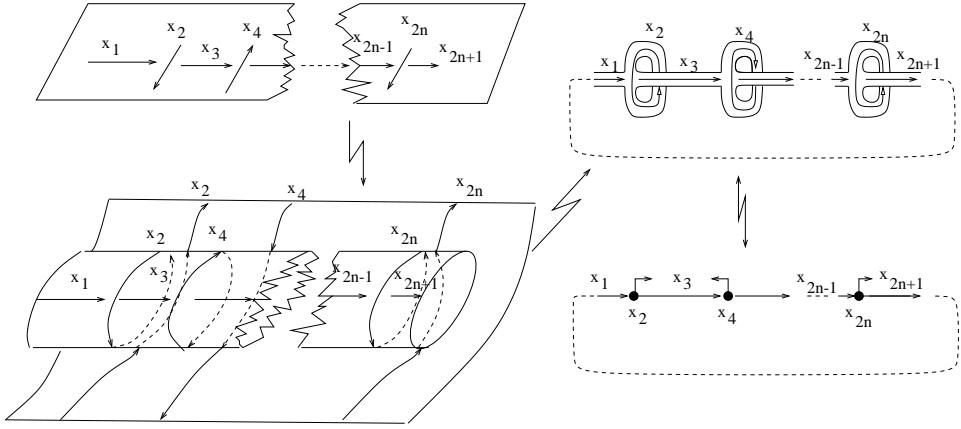}
}
\end{center}
\caption{The correspondence between shadow colored $0$-diagrams and
abstract colored $1$-diagrams}
\label{daisychain}
\end{figure}

\begin{sect}{\bf Extended Remark.\/}
\label{extended}
{\rm 
Let ${\mbox{\rm SC}}(1)$ denote the collection of shadow colored 
abstract  $1$-knot
diagrams, 
and let $C(2)$ denote the collection of
colored abstract
closed $2$-knot 
diagrams.
We have maps 
${\mathcal D}: C(2) \rightarrow {\mbox{\rm SC}}(1)$ and 
${\mathcal I}:{\mbox{\rm SC}}(1) \rightarrow C(2)$ that are
 defined as follows (see also Fig.~\ref{ab1eq2} for a local description).

Consider a colored 
closed 
abstract
$2$-knot diagram, $K^2= [f:M^2\rightarrow N^3]$. 
The map ${\mathcal D}$  assigns to 
$K^2$ 
the following shadow colored $1$-knot diagram.
The abstract 
$1$-knot  
diagram is 
the lower decker set in 
a  
regular neighborhood of
the lower 
decker set in the surface $M$ 
--- the 
lower decker points form the $1$-crossings and $2$-crossings.
Here, {\it lower decker set} 
(see \cite{CS:book}) means the preimage
of the double point set of the surface 
$M$ 
that are 
in the under sheet
and represented as the 
broken 
sheet 
in the abstract $2$-knot diagram. 
The $2$-crossing 
points of ${\mathcal D}(K^2)$
correspond to the triple points of the 
abstract $2$-knot diagram, 
 and the endpoint diagrams correspond to 
branch points (see \cite{CS:book} for details).  
We color  the arcs of ${\mathcal D}(K^2)$
 with the 
quandle elements that appear on the regions that contain the 
corresponding upper sheets. 
At a 2-crossing 
of ${\mathcal D}(K^2)$, 
the over-arc is colored by the 
color on the upper sheet, 
the under arcs are colored by the 
colors on the two portions of the corresponding middle sheet.
The 2-dimensional 
regions ${\mathcal D}(K^2)$ are colored by 
the quandle elements that are 
on the pieces of the surface that contain the arc diagram.

On the other other hand, 
the map ${\mathcal I}$ assigns to an abstract  
shadow 
colored 
1-knot 
diagram, $K^1$,
 a colored abstract $2$-knot diagram 
as follows.
Decompose the abstract $1$-knot diagram into 
crossing and endpoint diagrams. 
Construct an abstract  
$2$-knot 
diagram  
using the local correspondence 
between 
$k$-crossing $1$-diagrams and $(k+1)$-crossing $2$-diagrams, 
and between 
endpoints and branch points 
that is depicted in Fig.~\ref{ab1eq2}. 
Colors on the surface are determined as in 
the figure. 
In this way, a surface with 
double points on the boundary is constructed. 
The colors on the double points on the 
boundary just above and below the surface of the arc diagram
agree. Thus these double points can be joined together 
as they were
in
 Fig.\ref{genh3r3susp}.

Given $K^1$, the $1$-knot diagram 
${\mathcal D}({\mathcal I}  (K^1)  ) $
 differs from $K^1$ in that in 
addition to $K^1$, it contains unknotted 
unlinked colored components that may not be found in $K^1$.

In particular, given a shadow colored 
classical knot diagram, we obtain 
a 
colored abstract $2$-knot diagram. 
If $K$ is a 
knot, the $2$-knot diagram  
consists of a sphere (in the plane of the 
knot diagram) and a torus of the form 
$K\times S^1$. The $2$-knot diagram 
 is in  $S^2 \times S^1$. We call this the {\it suspension} 
of a classical knot. See Fig.~\ref{susp}.

A similar construction applies to give maps
${\mathcal D}: C(1) \rightarrow {\mbox{\rm SC}}(0)$ and 
${\mathcal I}:{\mbox{\rm SC}}(0) \rightarrow C(1)$ 
between the sets, $C(1)$,  
of colored abstract 
closed 
$1$-knot diagrams 
and ${\mbox{\rm SC}}(0)$ of 
shadow colored $0$-knot diagrams.
 The essense of this construction is found in \cite{Greene}
 (see also Fig.~\ref{daisychain}).

} \end{sect}

\section{Representing Homology Classes} \label{repsec}

 In Greene \cite{Greene} (see also \cite{FRS2,Flower}), rack 
cycles are represented by colored knot diagrams.
We generalize 
this  
method to quandle cycles, using 
end/branch point and hem diagrams. 
Specifically, for a fixed quandle 
$2$- or $3$-cycle we construct a 
colored or shadow colored 
diagram that represents the cycle. 
Moreover, if a cycle is a boundary, 
we construct an 
abstract $2$-knot 
diagram in a 3-manifold, 
the boundary of which is the given 
abstract $1$-knot 
diagram that represents the cycle. 
Care will be taken at endpoint diagrams 
to make this precise.
Finally, we represent $4$-cycles 
by shadow colored 
abstract 
$2$-knot diagrams
 with hems. We do not pursue boundaries in this case.

\begin{figure}
\begin{center}
\mbox{
\epsfxsize=4in
\epsfbox{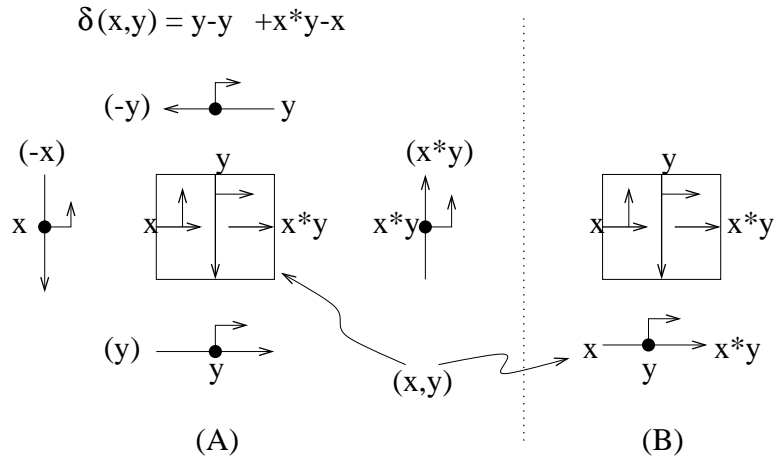}
}
\end{center}
\caption{A geometric representation of a generating $2$-chain and its boundary}
\label{2chains}
\end{figure}

\begin{figure}
\begin{center}
\mbox{
\epsfxsize=3.2in
\epsfbox{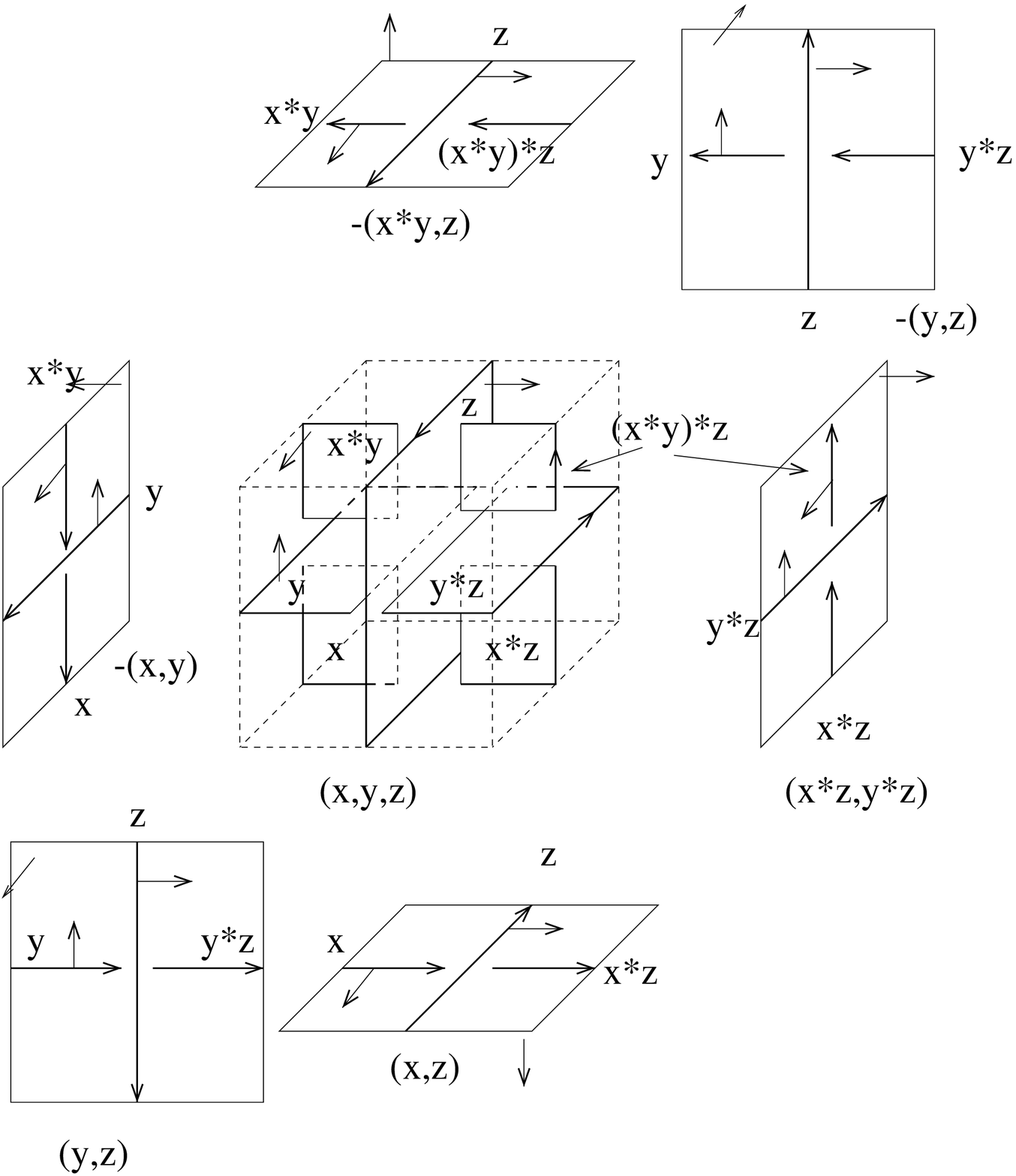}
}
\end{center}
\caption{A geometric representation of a generating $3$-chain and its boundary}
\label{3chains2}
\end{figure}

\begin{figure}
\begin{center}
\mbox{
\epsfxsize=2.8in
\epsfbox{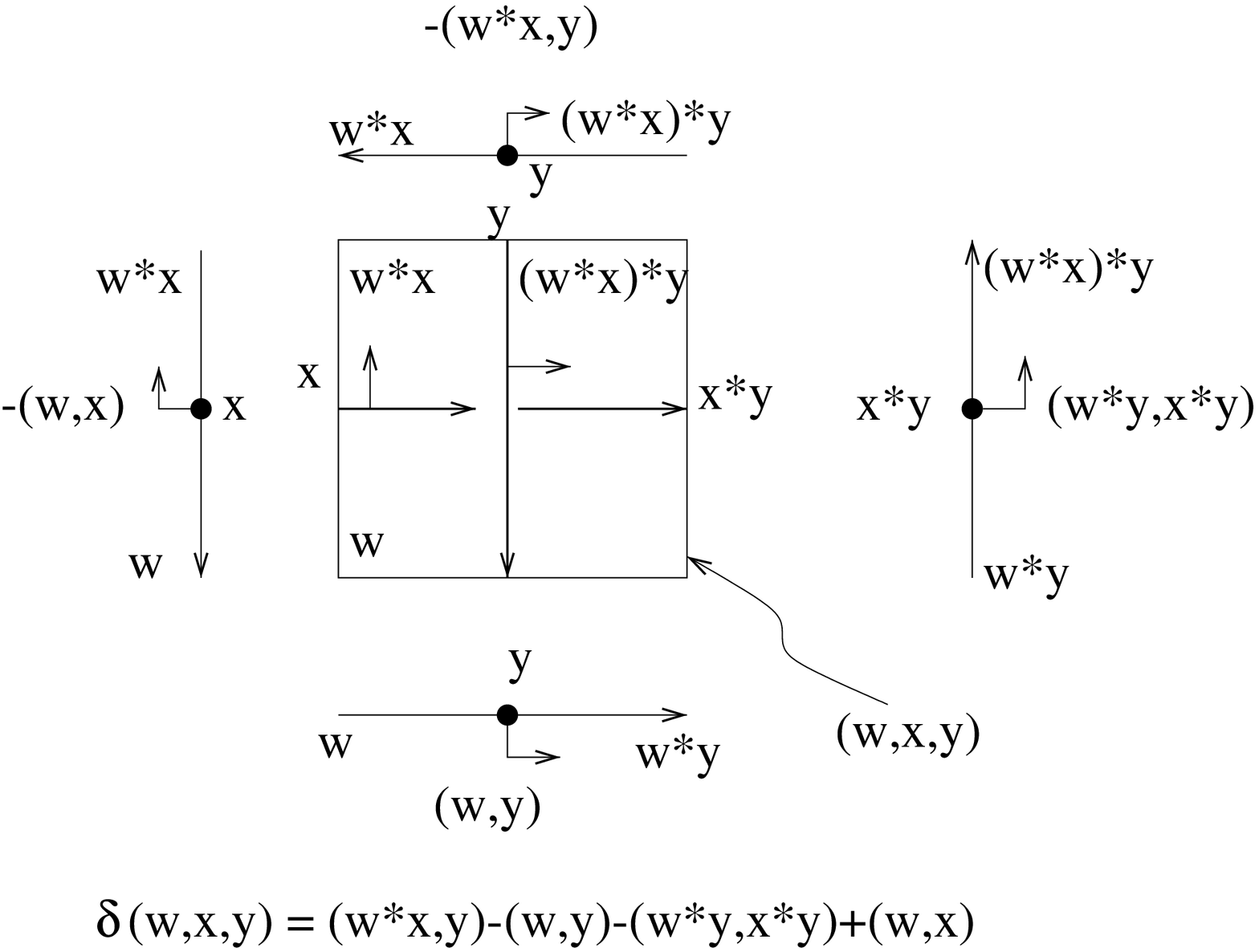}
}
\end{center}
\caption{A shadow representation of a generating $3$-chain and its boundary}
\label{3chain}
\end{figure}

\begin{sect}{\bf Lemma.} \label{tuple}  
(1) There is a one-to-one correspondence between
the set of 
$(n+1)$-tuples of quandle elements and 
the set of quandle colored  positive $(n+1)$-crossing
$n$-diagrams.

(2) There is a one-to-one correspondence between 
the set of $(n+1)$-tuples of 
quandle elements and shadow colored positive 
$n$-crossing $n$-diagrams.
\end{sect}
 {\it Proof.} The correspondence (1) for $n=1$ is   
illustrated in  
 Fig.~\ref{2chains} (A). 
The correspondence
(2) for $n=1$ is illustrated 
in 
Fig.~\ref{2chains} (B), 
where the correspondence between the shadow colors and the lower decker
set is illustrated. 
The correspondence (1) for $n=2$  is illustrated in the center of Fig.~\ref{3chains2}. 
The correspondence (2) for $n=2$ is illustrated in Fig~\ref{3chain}.

In general, let $X$ be a quandle and
consider an $(n+1)$-tuple $(x_1, \ldots , x_{n+1}) \in X^{n+1}$. 
Recall that 
$E_j$ is the coordinate $n$-disk in $[-1,1]^{n+1}$ 
whose $j\/$th coordinate is $0$. Color the region of $E_j$ at which
the $\ell$\/th  coordinates are less than $-\varepsilon$ for $\ell <j$
 with $x_{n+2-j}.$ Color the remaining regions of $E_j$ in such a way that
the quandle condition holds at each double point of the diagram. 

Such a coloring can be uniquely
extended to a shadow coloring by coloring the 
$(n+1)$-region 
of $E$ at which all the coordinates are less than $-\varepsilon$ by 
$x_0$. Thus associated to an $(n+1)$-tuple there is a colored $(n+1)$-crossing diagram or a shadow colored $n$-diagram. 

By choosing the color $x_{n+2-j}$  
from the region of $E_j$  at which the $\ell\/$th coordinates are 
less than $- \varepsilon$ 
 for $\ell < j$, we construct an $(n+1)$-tuple.
A shadow colored diagram works similarly. $\Box$

\begin{sect}{\bf Scholium.\/} The boundary 
$\partial (x_1, \ldots, x_{n+1})$ 
of a generating rack chain
$(x_1, \ldots, x_{n+1})$ corresponds 
to the boundary of the quandle colored
$(n+1)$-crossing diagram. 
\end{sect}
 {\it Proof.} The correspondence is illustrated in the 
Figs.~\ref{2chains} (A) 
 and \ref{3chains2}. We leave the rest of the details to the reader.
$\Box$

\begin{sect}{\bf Remark.\/}
\label{itsacycle}
{\rm 
Colored  and  shadow colored 
abstract knot diagrams represent chains ---
formal-sums of 
signed 
$(n+1)$-tuples 
(or $(n+2)$-tuples) of quandle elements 
corresponding to the $0$-dimensional 
crossings as in Lemma~\ref{tuple}. 
Recall that the $0$-dimensional crossings of an
abstract $n$-knot diagram are the
$(n+1)$-crossing points. 
Each such (shadow) colored crossing point  in the diagram
has an associated sign.
The chain determined by the diagram is the sum of 
these  signed 
$(n+1)$-tuples (or $(n+2)$-tuples) 
taken over all the $0$-dimensional crossings.
This is the {\it chain represented by a (shadow) colored diagram}.

An abstract diagram (of a closed manifold $M$)
that has no exceptional points has its
$1$-dimensional crossing set closed. Therefore  
the sum of the representative 
chain is a 
rack 
cycle in this case \cite{FRS1}.

In low dimensions, 
shadow colored exceptional points or  colored branch points
represent degenerate chains. 
Specifically, a shadow colored endpoint or a colored 
branch point represents a chain of the form $\pm (a,a)$, a shadow colored 
branch point represents a chain of the form $\pm (a,b,b)$
and a shadow colored hem $2$-crossing point represents a chain
$\pm (a,a,b)$. The signs are determined by the direction of the arrow
along the arc that originates or 
terminates at the exceptional point. 

Thus 
a shadow colored abstract $0$-knot diagram represents a 
quandle $2$-cycle (for example, Fig.~\ref{example0}), as does a
colored closed abstract $1$-knot diagram (for example, Fig.~\ref{ssvsv1}).
A shadow colored  
abstract $1$-knot diagram represents a 
{\it quandle} $3$-cycle as does a colored
closed abstract 
$2$-knot diagram (for example, Fig.~\ref{genh3r3susp}). 
And a shadow colored abstract $2$-knot diagram
represents a {\it quandle} $4$-cycle. 
For example the 2-twist-spun trefoil
may be shadow colored by $R_3$ (see \cite{CJKLS}, and apply Lemma~\ref{genh3r3}).
This is the context of Theorem~\ref{representation} (1a) and (1b).
}\end{sect}

\begin{sect} {\bf Examples.\/}
{\rm 
The  $R_3$
$2$-chain $(\alpha,\beta)+(\gamma, \alpha)+(\beta,\gamma),$
represented by the diagram in Fig.~\ref{example0} is a cycle
in $Z^{\rm Q}_2(R_3)$.

The $R_3$ $2$-chain
represented by 
the 
colored abstract $1$-knot diagram in Fig.~\ref{ssvsv1} is the cycle
$(\alpha,\beta)+(\beta,\gamma)-(\beta,\alpha) \in Z_2^Q(R_3;{\bf Z})$.

It is known \cite{CJKLS} that $H_2^Q(R_3;{\bf Z})=0$, 
so that the above 
examples are in fact boundaries.

The $3$-chain
represented by 
the 
colored abstract $2$-knot diagram in Fig.~\ref{genh3r3susp} is
the cycle 
$(\alpha,\beta,\gamma)+(\alpha,\gamma,\alpha) \in Z_3^Q(R_3;{\bf Z})$.
The corresponding 
(in the sense of \ref{extended})
shadow colored abstract $1$-knot diagram is 
depicted at the right bottom of the figure, which 
is the same as Fig.~\ref{genh3r3}. 
It is known \cite{CJKS1} 
that  $H_3^Q(R_3;{\bf Z}) \cong {\bf Z}_3$, and the above $3$-cycle 
is a generator.

} \end{sect}

\begin{sect}{\bf Theorem.\/}
\label{representation}
Let $X$ denote a quandle. 

(1a) Let $n=1,2$.
Any colored 
closed 
abstract $n$-knot diagram 
represents a quandle $(n+1)$-cycle in $Z^{\rm Q}_{n+1}(X; {\bf Z})$. 

(1b) Let $n=1,2,3.$
Any  shadow colored abstract $(n-1)$-knot 
diagram (possibly 
with exceptional points if $n=2,3$) 
represents a quandle 
$(n+1)$-cycle in $Z^{\rm Q}_{n+2}(X;{\bf Z})$. 
 
(2a)
Let $n=1,2$. 
 Let $\eta \in Z^{\rm Q}_{n+1}(X; {\bf Z})$.  
Then there is a colored $n$-knot diagram, 
${\mathcal D}^n_\eta$, that 
represents $\eta$.

(2b) Let $n=1,2,3.$
 Let $\eta \in Z^{\rm Q}_{n+1}(X; {\bf Z})$.  
Then there is a 
shadow
colored  
$(n-1)$-knot
diagram
(possibly 
with exceptional points for  $n=2,3$), 
${\mathcal SD}^{n-1}_\eta$,
that represents $\eta$.
\end{sect}
{\it Proof.\/}
Statements (1a) and (1b) (for $n=2,3$), 
follow from Remark~\ref{itsacycle}.
For $n=1$ statement (1b) is also easy: A shadow colored $0$-knot diagram represents a rack $2$-cycle which is also a quandle $2$-cycle.

Consider statement (1b) for $n=2$, 
the endpoints of shadow arc diagrams
represent chains of the form 
$(a,a)$ which are trivial in quandle 
homology. 
Thus shadow colored arc 
diagrams represent quandle $2$-cycles.

Branch points
of shadow colored abstract $2$-knot 
diagrams represent chains of the form
$(a,b,b)$ while hem $2$-crossing points 
represent chains of the form $(a,a,b)$. 
Both such chains are trivial in quandle homology.
This proves (1b) for $n=3$.

Consider a quandle 
$2$-cycle $\eta=\sum_j \epsilon_j
\vec{x}_j \in Z_2^Q(X; {\bf Z})$ 
where $\vec{x}_j= (x_1^j,x_2^j)$  
and $\epsilon_j  = \pm 1$.
For each  $(x_1^j,x_2^j)$, 
construct a colored $2$-crossing $1$-diagram as in Lemma~\ref{tuple}
(1); the sign of the crossing is $\epsilon_j$.
 The boundary
of such a chain is 
$\partial (x_1^j,x_2^j)= 
 x_1^j - x_1^j  * x_2^j $. 
 Since $\eta$
is a cycle the sum of these boundary 
terms adds to $0$. Thus we can 
interconnect the under-arcs to 
form a collection of simple closed curves
  of under-arcs. The over-arcs each 
can be joined end-to-end locally 
to form the 
abstract diagram
as depicted in Fig.~\ref{daisychain} right top.
 For example,
 Fig.~\ref{ssvsv1} illustrates an abstract 
diagram that represents
the cycle $( \alpha, \beta) + (\beta, \gamma) - (\beta, \alpha)\in Z^{\rm Q}(R_3)$. 
By Remark~\ref{extended}, a desired shadow colored $0$-knot diagram 
for (2b) is obtained from the above constructed 
 representative $1$-knot diagram.
Thus statements (2a) and (2b) hold for $n=1$.
Note that the way canceling pairs are joind together is not unique.

Now suppose that $\eta = \sum_j \epsilon_j \vec{x}_j 
\in Z_3^Q(X; {\bf Z})$  
(where
$\vec{x}_j= (x_0^j,x_1^j,x_2^j)$ 
and $\epsilon_j = \pm 1$) 
is a $3$-cycle. 
We represent each $\epsilon_j \vec{x}_j$ by a shadow colored
$2$-crossing  
$1$-diagram 
as in Lemma~\ref{tuple} (2); the sign of the crossing is
$\epsilon_j$.
 The boundary of such a chain is the sum of 
four $2$-chains. Each $2$-chain is 
represented as a 
shadow colored $0$-knot diagram. Take the cartesean product
of such a diagram with the unit interval to form a shadow colored 
$1$-crossing $1$-knot diagram. 
Such a diagram is  a $1$-handle
the attaching region of which is 
identified with the segments 
on the boundary of the squares
that represent the chains $\vec{x}_j$. 
Thus we can attach $1$-handles to these squares.
Since $\eta$ is a quandle cycle, the boundary of squares that are
not attached to these $1$-handles 
have the same color assigned to 
the arc and the region.
Cap  
each of  these boundary segments  by a
colored endpoint diagram. This constructs a colored
abstract $1$-knot diagram 
with endpoints 
representing $\eta$.  
 To construct an abstract $2$-knot diagram 
that represents $\eta$, we 
follow the procedure  
outlined in 
Remark~\ref{extended}. 
This gives statement (2a) and (2b) for $n=2$. 

\begin{figure}
\begin{center}
\mbox{
\epsfxsize=3in
\epsfbox{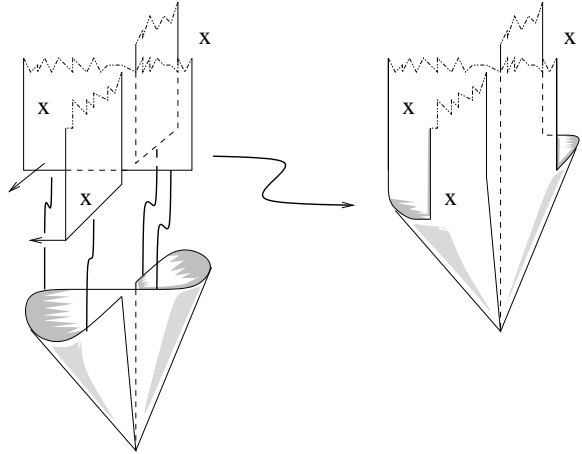}
}
\end{center}
\caption{Capping off with a branch point} 
\label{duncehat}
\end{figure}

If $\eta = \sum_j \epsilon_j \vec{x}_j$  (where
$\vec{x}_j= (x_0^j,x_1^j,x_2^j,x_3^j)$ 
and $\epsilon_j = \pm 1$) 
is a $4$-cycle, then we construct 
a shadow colored 
abstract $2$-knot diagram that may have hems  
following a procedure analogous 
to the ones above. Specifically,
for each chain $\vec{x}_j$ we 
construct a shadow colored 
triple point diagram where the 
colors on the regions away from which normals point is $x_i^j$.
The boundary of the chain $\vec{x}_j$ 
corresponds to the shadow colored $2$-crossing diagrams on the 
boundary of the colored cubes. 
We attach shadow colored $1$-handles 
between pairs  of square faces to cancel like terms.
Since $\eta$ 
is a quandle cycle, its boundary consists of terms of the form
$(a,a,b)$ and $(a,b,b)$.
 There are faces that remain  
representing 
such terms. 
In case the term is of the form $(a,b,b)$,
 we attach a branch point as in Fig.~\ref{duncehat}. 
In case the term is the form $(a,a,b)$, then we extend the over-sheet
 to create a hem diagram. The double point of the hem diagram 
represents 
 this chain.
This completes the proof of 
(2b) for $n=3$.
$\Box$

\begin{figure}
\begin{center}
\mbox{
\epsfxsize=4in
\epsfbox{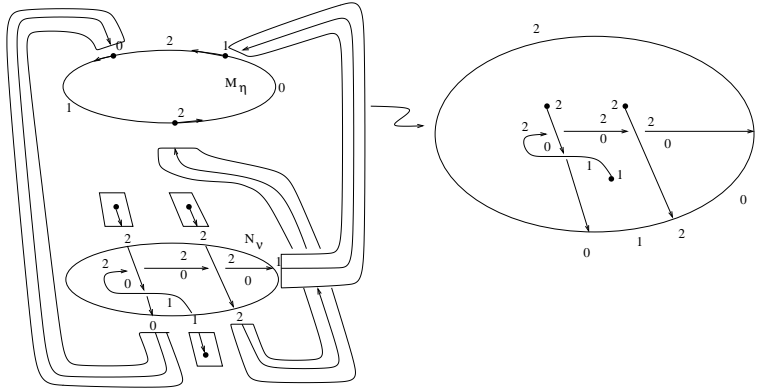}
}
\end{center}
\caption{Cobordism for shadow colored $0$-diagrams}
\label{sh1}
\end{figure}

\begin{sect}{\bf Remark.\/} 
\label{actualb}
{\rm 
We included condition (1) 
of Definition~\ref{absdef} 
for expository convenience. 
Otherwise boundary points might be 
confused with hems or endpoints.
In the following theorem, 
we 
use 
{\it abstract diagrams with actual boundaries}
 --- diagrams for which condition (2) is satisfied, 
but not necessarily condition (1).
The notions of colored and shadow colored diagrams
extend to diagrams with actual 
boundary in a straight-forward way.

In ordinary homology theory, boundaries of complexes correspond
to boundary terms  of chains. This is our motivation of 
extending the abstract diagrams
to those with boundaries, as will be seen in the following theorem.
The exceptional points that are defined above are not regarded as
boundaries in 
quandle homology,
as they were 
introduced to represent degenerate chains, instead of 
boundary terms.

For example, an 
abstract $1$-knot diagram 
with actual boundary 
$[ f: M \rightarrow N]$ has
two types of endpoints. One is the endpoint diagram defined already.
The {\it actual boundary end point}
lies on the boundary $\partial N$. 
When a diagram is shadow colored by a quandle, the former 
(the endpoint diagram) corresponds to a 
degenerate chain, 
and the latter
corresponds to the 
sum of the 
boundary terms in $\partial (x_0, x_1, x_2)$. 
See 
the right-side of Fig.~\ref{sh1}, for example.

In case a shadow  
colored 
$2$-knot 
diagram 
with actual boundary
has a hem, 
the hem may have a boundary point in the 
form of an endpoint diagram, but 
the hem 
itself is not  
regarded as  
the actual boundary.

} \end{sect}

\begin{sect}
{\bf Theorem.\/} 
 Let $n=1,2$.  
 If the given cycle,  $\eta \in Z^{\rm Q}_{n+1}(X; {\bf Z})$, 
is a boundary, 
so $\eta= \partial \nu$, 
for $\nu  \in C^{\rm Q}_{n+2}(X; {\bf Z})$,   
then 
${\mathcal D}^n_\eta$ (or ${\mathcal S D}^{n-1}_\eta$)
is the boundary of a colored 
$(n+1)$-knot 
diagram, ${\mathcal D}^{n+1}_\nu$, 
(or shadow colored $n$-knot diagram, ${\mathcal SD}^{n}_\nu$,) 
with actual boundary 
(and with hems when appropriate)
 that represents $\nu$. 
\end{sect}
{\it Proof.} 
Suppose that the  quandle $(n+1)$-cycle $\eta = \partial \nu$. 
In general, there is a
colored abstract $(n+1)$-knot diagram (or shadow colored 
abstract $n$-knot diagram) with actual boundary 
that represents $\nu$.
For example,
 take the 
disjoint union 
of 
colored $(n+2)$-crossing $(n+1)$-diagrams 
(or shadow colored $(n+1)$-crossing $n$-diagrams)
that represent the chains which constitute $\nu$.
Let $({\mathcal S}){\mathcal D}_\nu $ denote such a diagram, where 
the parenthesis represents that it is either a colored diagram 
(without ${\mathcal S}$) or
 shadow colored diagram (with ${\mathcal S}$).
We also use the  notation
$({\mathcal S}){\mathcal D}_\nu
= \left[ g: M_\nu \rightarrow N_\nu \right]_{ {\cal C}_\nu } $ 
to specify the map $g$.
Similarly, let $({\mathcal S}){\mathcal D}_\eta 
= \left[f: M_\eta \rightarrow N_\eta \right]_{ {\cal C}_\eta }$
denote a given 
colored diagram that represents $\eta$.

The sum of the (shadow) colored crossings
of $\partial ({\mathcal S}){\mathcal D}_{\nu} $ 
may differ from the sum of the 
(shadow) colored crossings on 
$({\mathcal S}){\mathcal D}_\eta$
in that either may include 
degenerate chains 
or canceling terms. 
That is, there can be crossings that represent
chains of the form
$(a,a)$ in case $n=1$, or chains of one of the forms 
$(a,a,b)$  or $(a,b,b)$ if $n=2$. 
Or there may be a canceling pair of  
(shadow) colored crossings.

Take the product $({\mathcal S}){\mathcal D}_\eta \times [0,1]$ 
which is $g(M_{\eta}) \times [0,1] \subset N_{\eta}  \times [0,1]$
with the given  coloring ($\times [0,1]$),
where the original $({\mathcal S}){\mathcal D}_\eta $ is regarded as 
embedded in $N_{\eta}  \times \{ 0 \}$. 
For each  pair of canceling (shadow) colored crossings, 
join them
by a (shadow) colored $1$-handle 
(with coloring induced from the crossings), and
cancel them.
Such   
$1$-handles 
are  
attached on   $N_{\eta}  \times \{ 1 \}$, 
and give 
a new diagram $[ g' : M_{\eta}' \rightarrow N_{\eta}']$
with a (shadow) color 
such that $\partial N_{\eta}'$ consists 
two pieces, 
one of which $\partial _1 N_{\eta}'$ contains the original  $g(M_{\eta})$, 
and the other  $\partial_2 N_{\eta}'$ 
contains  $g(M_{\eta})$ with all canceling pairs eliminated.
Perform the same process on $\partial N_{\nu}$  
(we continue to use the same notation $N_{\nu}$ after this process),
so that we now assume that  $\partial _2 N_{\eta}'$ and 
 $\partial N_{\nu}$ do not contain canceling pairs.
Next we cancel degenerate chains, case by case.

Case 1.
Suppose that $\eta$ is a $2$-cycle
represented by a shadow colored  $0$-knot diagram, 
${\mathcal SD}_\eta=[f:M_\eta \rightarrow N_\eta]_{ {\cal C}_\eta } $.  
Then $M_\eta$ is a $0$-dimensional 
manifold, $N_\eta$ is a $1$-dimensional 
manifold. We call the components of $M_\eta$ oriented {\it vertices}. 
A $2$-chain $(a,b)$ is represented by a 
shadow colored $1$-crossing $0$-diagram 
 where 
$b$ is the color of the vertex 
and $a$ is the color on the arc 
away from which the normal arrow to the 
vertex points. 
The description of 
$\partial{\mathcal SD}_\nu$ 
is the same.
For each degenerate vertex (shadow colored as $(a,a)$), 
we construct a shadow colored endpoint diagram with color 
$a$ on the $2$-dimensional 
region and color $a$ on the edge.
The orientation of the arc is determined by the sign of the 
$0$-crossing. Each such endpoint diagram is a $0$-handle that
is disjoint from $N_\eta \times [0,1]$ and $N_\nu$.
Now attach a $1$-handle (in the guise of a shadow colored
$1$-crossing $1$-diagram)
 between the degenerate crossing and the endpoint diagram. 
In this way, we may assume that there are no degenerate $2$-chains
among the $1$-crossings representing the 
summands of $\eta$ and  $\partial \nu$.
See Fig.~\ref{sh1} for an example.

Case 2.
Suppose 
that $\eta$ is a $2$-cycle that is represented 
by a colored $1$-diagram.
Each degenerate colored $2$-crossing 
in 
${\mathcal D}_\eta$ or
$\partial {\mathcal D}_\nu$ is a 
crossing at which all the colors are 
the same; thus a degenerate crossing represents a chain 
for the form $(a,a)$. 
For each such crossing we construct a colored 
branch point diagram that 
functions as a $0$-handle. 
The color on the surface is $a$.
This handle 
is attached  to the diagram via a 
$1$-handle in the guise of a colored $2$-crossing
$2$-diagram. Therefore, as above we can assume 
that there are no degenerate
$2$-chains among the crossings representing
the summands $\eta$ and $\partial \nu$.
See Fig.~\ref{duncehat}.

Case 3.
Suppose 
$\eta$ is a $3$-cycle that is represented by a shadow colored 
$1$-knot diagram, ${\mathcal SD}_\eta$.   
The degenerate chains are of the form $(a,a,b)$
and $(a,b,b)$. 
The former are represented by colored crossings with color $a$ on one face,
color $a$ on one of the lower arcs, and color $b$ on the upper arc. The latter
are represented by shadow colored crossings where $a$ is the color on one face, and  $b$ is the color on all edges. 
For each degenerate chain of the form 
$(a,a,b)$, we insert a shadow 
colored hem $2$-crossing diagram as a $0$-handle 
disjoint from $N_{\eta}'$ 
 and $N_\nu$. For each degenerate
chain of the form 
$(a,b,b)$
we insert a shadow colored
branch point diagram as a $0$-handle. Then the $0$-handles 
are attached to either  $\partial _1 N_{\eta}'$ 
 or $\partial N_\nu$ via 
$1$-handles that are $2$-crossing $2$-diagrams. 
The $1$-handles are attached to the $0$-handles at the crossing diagram on the boundary.
Thus in this case, we may assume that there are no degenerate 
$3$-chains
among the crossings 
representing the summands of  $\eta$ and $\partial \nu$.

Next, we consider the cases when the colored 
crossings on the resulting manifolds 
(that we again denote by $\partial _1 N_{\eta}'$ 
 and $\partial N_\nu$)
are in one-to-one correspondence. In each case, we can attach 
$1$-handles between pairs of similarly (shadow) 
colored crossings. The $1$-handles are of the form:
(a) shadow colored $n$-crossing  $n$-diagrams 
in case 1 and 3;  
(b) colored $(n+1)$-crossing 
$(n+1)$-diagrams in case 2.

After all of the crossings that represent summands of $\eta$ have been 
glued to the chains representing $\nu$, we attach (shadow) colored 
$2$-handles in the guise of 

\noindent
(1) shadow colored $0$-crossing $0$-diagrams 
 in case 1.

\noindent
(2)  $0$-crossing $1$-diagrams 
in case 2.

\noindent 
(3) shadow colored $1$-crossing 
$2$-diagrams 
in case 3.
 
This completes the proof. $\Box$

\section{Equivalence of 
Colored 
and Shadow Colored Diagrams} \label{equivsec}

So far we have introduced abstract diagrams, 
colorings, and shadow colorings 
thereof. Here we discuss moves to such diagrams. 

\begin{figure}
\begin{center}
\mbox{
\epsfxsize=5in
\epsfbox{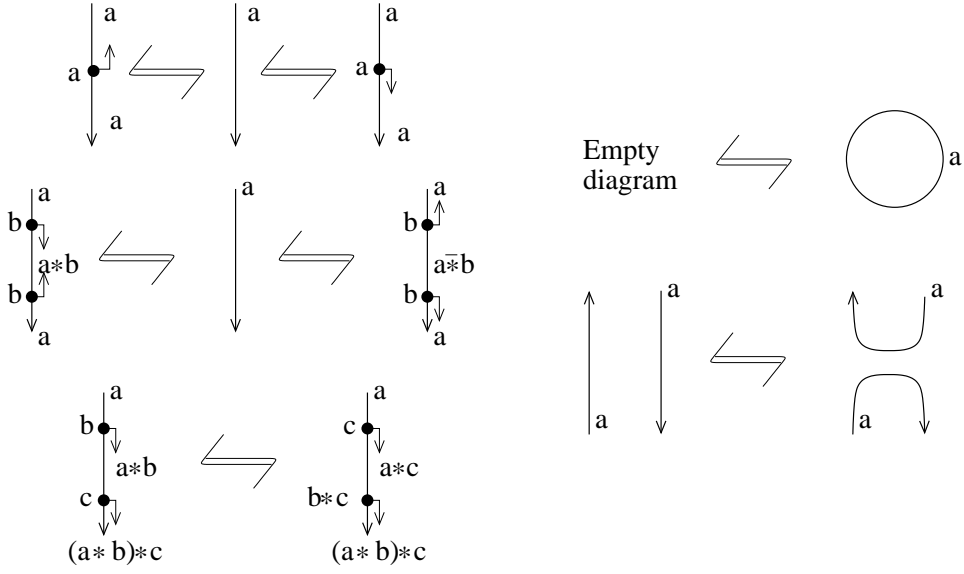}
}
\end{center}
\caption{Moves for shadow colored 0-diagrams}
\label{0dcobmoves}
\end{figure}

\begin{sect} {\bf Theorem.\/} 
For $i=0,1$, 
let ${\mathcal SD}_i$ 
 denote 
shadow colored abstract $0$-knot 
diagrams where the color set is the quandle 
$X$. Suppose that
$[{\mathcal SD}_0]=[{\mathcal SD}_1] \in H^{\rm Q}_{2}(X; {\bf Z})$.
Then ${\mathcal SD}_0$ 
 can be 
obtained from ${\mathcal SD}_1$ 
by a finite sequence of moves 
taken from those depicted in Fig.~\ref{0dcobmoves}.
\end{sect} 
 {\it Proof.\/} 
Let ${\mathcal SD}_i $ be represented by maps with colorings
$[f_i: M_i \rightarrow N_i]_{{\cal C}_i}$ as in the proof of the 
preceding theorem.
By Theorem~\ref{representation},  
$[f_0:M_0 \rightarrow N_0]_{{\cal C}_0}$  and   $[f_1:M_1 \rightarrow N_1]_{{\cal C}_1}$ 
cobound a shadow 
colored abstract $1$-knot diagram $M$
in a $2$-manifold $N$ 
with actual boundary. 
Let $F: N \rightarrow [0,1]$ be a smooth function 
such that 
$F^{-1}(i) = N_i,$  for $i=0,1.$
We may assume (after a small perturbation if necessary) 
that $F$ satisfies the following conditions:

\noindent
(1) $F$ is transverse at $0$ and $1$.

\noindent
(2) $F$ is generic on $M$, $N$, $\partial M$, $\partial N$, 
and on all the self intersections and singularities of $M$.

Thus $F$ has 
isolated 
Morse critical points on all the sets listed in (2),
at distinct critical values. 
By taking the inverse images $F^{-1}(h -\epsilon)$ and 
$F^{-1}(h + \epsilon)$ at every critical value $h$,
we obtain a sequence of moves. Hence we obtain the result 
by classifying these Morse critical points as follows.

They are classified into two categories: Reidemeister moves
for $0$-dimensional knot diagrams as singular sets and critical points of $M$
(Fig.~\ref{0dcobmoves} left $3$ figures), and Morse critical points of $N$ 
(Fig.~\ref{0dcobmoves} right $2$ figures).
The 
left 
$3$ figures correspond to the endpoints of $M$, 
maxima/minima of $M$, and transverse double points of $M$, respectively.
The 
right 
$2$ figures correspond to the maxima/minima of $N$,
and saddle points of $N$, respectively. These exhaust generic critical
points of $f$, and the theorem follows. $\Box$

\begin{figure}
\begin{center}
\mbox{
\epsfxsize=3in
\epsfbox{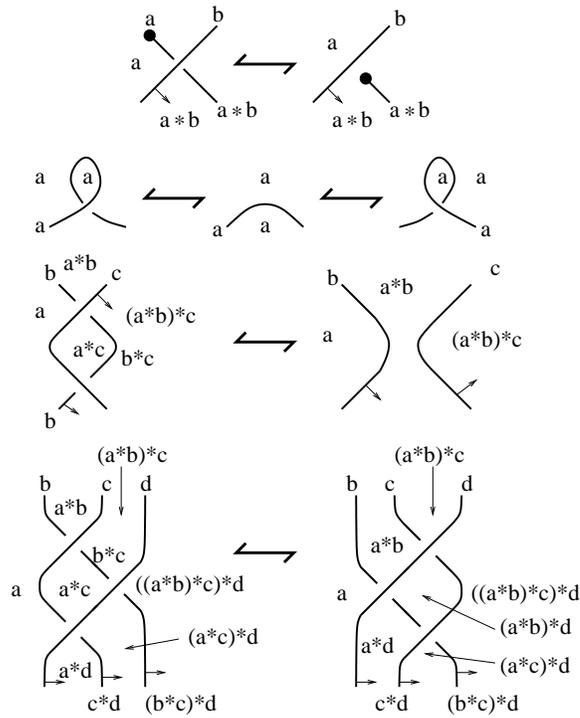}
}
\end{center}
\caption{Moves for shadow colored diagrams, Part I (Reidemeister moves)}
\label{shmoves2}
\end{figure}

\begin{figure}
\begin{center}
\mbox{
\epsfxsize=4.5in
\epsfbox{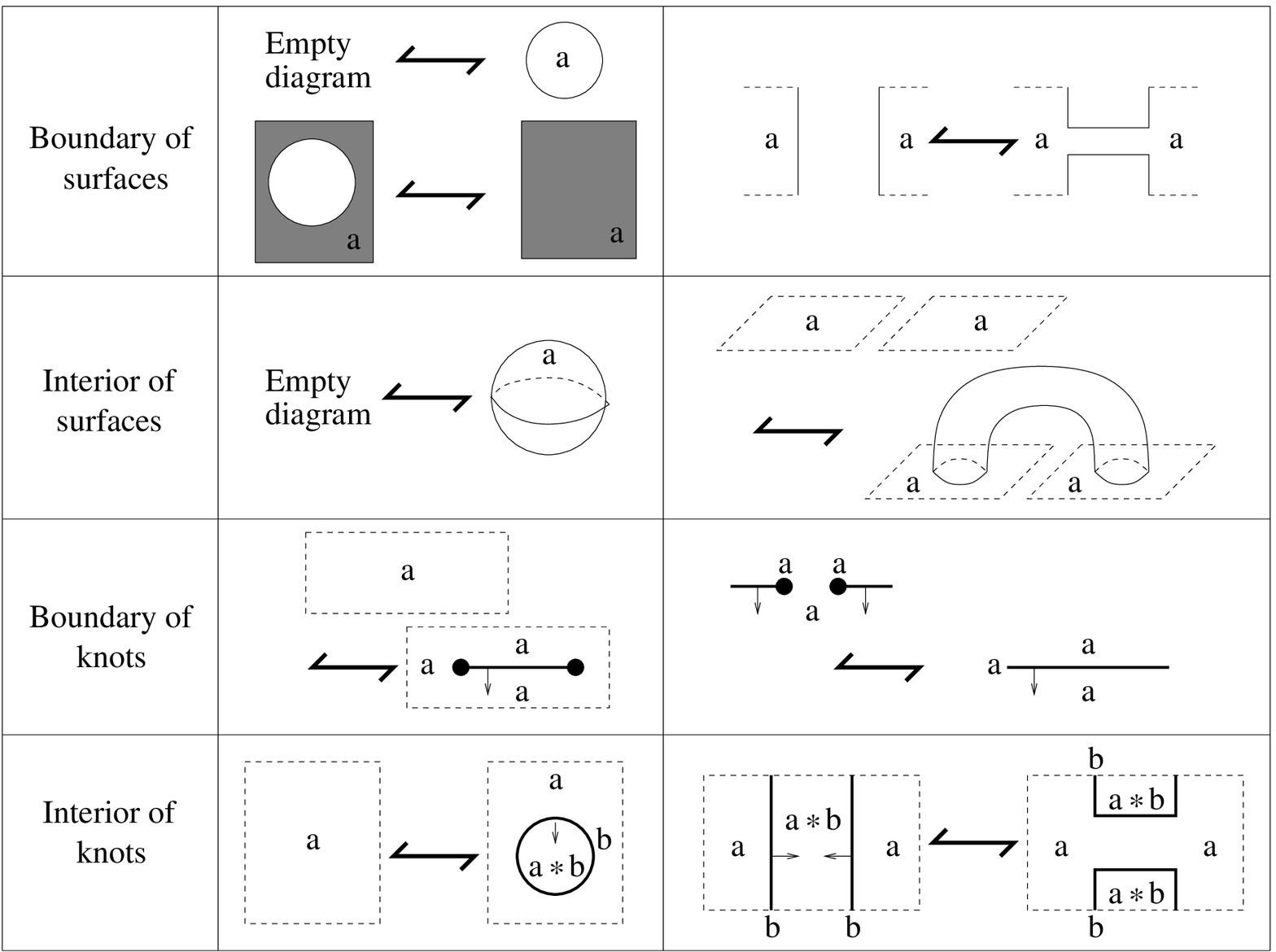}
}
\end{center}
\caption{Moves for shadow colored diagrams, Part II (Morse critical points) }
\label{shmoves1}
\end{figure}

\begin{sect} {\bf Theorem.\/} 
\label{morse2}
Let $X$ be a quandle.
For $i=0,1$, 
let  ${\mathcal SD}_i=[f_i:M_i \rightarrow N_i]_{ {\cal C}_i}$ 
denote shadow colored $1$-knot diagrams  that represent
the same homology class in $H^{\rm Q}_3(X;{\bf Z})$. 
Then $[f_1:M_1 \rightarrow N_1]_{ {\cal C}_1}$ can be obtained from 
$[f_0:M_0 \rightarrow N_0]_{ {\cal C}_0}$ by a finite sequence of moves  taken from those 
depicted  in 
Figs.~\ref{shmoves2} and \ref{shmoves1}.
\end{sect} 
{\it Proof.\/}
By Theorem~\ref{representation},  
$[f_0:M_0 \rightarrow N_0]_{ {\cal C}_0}$  and   $[f:M_1 \rightarrow N_1]_{ {\cal C}_1}$
cobound a 
shadow
colored abstract $2$-knot diagram $M$
in a $3$-manifold $N$ with actual boundary.  
Let $F: N \rightarrow [0,1]$ be a smooth function 
such that 
$F^{-1}(i) = N_i$, for $i=0,1$.
We may assume (after a small perturbation if necessary) 
that $F$ satisfies the following conditions.

(1) $F$ is transverse at $0$ and $1$.

(2) $F$ is generic on $M$, $N$, $\partial M$, $\partial N$, 
and on all the self intersections and singularities of $M$.

Thus $F$ has Morse 
isolated
critical points on all the sets listed in (2),
at distinct critical values.
The proof proceeds as follows in a similar way as in the preceding theorem. 

 The singularities and critical points of  $ M$
are listed in  Fig.~\ref{shmoves2} and are Reidemeister moves.
{}From top to bottom they represent the intersection of the boundary points 
and interior of $M$, branch points, maxima/minima of double curves of $M$,
and triple points, respectively.

The Morse critical points as handle moves are 
listed in Fig.~\ref{shmoves1}.
{}From top to bottom, they are critical points of 
$\partial N$, Int$N$, $\partial M$, and Int$M$.
The critical points of $\partial N$ are maxima/minima (the left entry)
or saddle points (the right). From the 
point of view of the boundary
 $1$-manifold, they correspond to handles of indices $0/2$ and $1$,
 respectively.
The critical points of Int$N$ are similar, and depicted in the second
 row 
left and right.
The critical points of  $\partial M$ are the 
creation or deletion
of a pair of points. 
There are two types,  left and right, of these $0/1$-handles,
in relation to interior points. 
The bottom entry 
in the figure
illustrates the 
critical points of the interior of $M$. 
The 
Theorem follows as these exhaust generic 
singularities 
and critical points.
$\Box$

\begin{sect}{\bf Scholium.\/}
Colored abstract  closed 
1-knot diagrams represent the same $2$-dimensional
quandle 
homology class if and only if 
one is obtained from the other 
by a finite number of moves 
that are taken from those depicted in Figs.~\ref{shmoves2} and \ref{shmoves1}
 where:

\noindent
{\em (1) those moves that involve 
endpoints 
are excluded; and }

\noindent
{\em (2) colors on the 2-dimensional regions are ignored.}
\end{sect}
 {\it Proof.} 
Imitate the proof of Theorem~\ref{morse2} 
replacing shadow colored with colored throughout. Since $M_i$ are closed, we may assume that the cobounded 
surface $M$ has no hems. Thus the 
moves involving endpoints are excluded. 
$\Box$

\begin{sect} {\bf Remark.\/} 
{\rm
For higher dimensions, 
the authors expect 
similar theorems.
The classifications of moves for boundaries, however, become
subtle. For one dimensional higher case, for example, 
the moves include Roseman moves (see for example \cite{CS:book})
with quandle colors, Morse critical points of the ambient space $N$
and $\partial N$, and singularities and critical points 
involving hems will provide the set of moves. 
The moves for hems are depicted in Fig.~\ref{funnystuff}.
} \end{sect}

\begin{sect}{\bf Discussion.\/}
{\rm 
Abstract $1$-knot diagrams that have no endpoints 
correspond to virtual knot diagrams (see \cite{KK}).
In \cite{KK}, a set of moves to abstract knot diagrams was given, and 
it was shown that up to these moves the set of abstract diagrams is 
equivalent to virtual knots up to virtual Reidemeister moves.
Also, the correspondence was made by thickening the virtual knots and
making them abstract knots. However, the moves to abstract knots are not to be 
confused with virtual Reidemeister  moves. 
Virtual Reidemeister moves for thickened virtual knots are not equivalence 
moves for abstract knots thus obtained. 
For example, type II moves to thickened virtual knots do not 
necessarily preserve 
shadow coloring.
We defined moves to colored and shadow colored diagrams to avoid 
diagrams that are type II equivalent but that have different quandles.
Therefore, quandle homology provides an
 invariant of colored diagrams under (shadow) colored Reidemeister moves
listed above.

}\end{sect}

\section{The Quandle of an Abstract Knot and Examples} \label{fundsec}

The (fundamental) 
quandle of an abstract $n$-knot 
diagram is defined in \cite{KK}. It is generated by the 
$n$-regions of the diagram; the relations in the quandle can be 
read  from the $2$-crossings. 
It generalizes the knot quandle of classical knots 
(\cite{Joyce,Matveev}). See \cite{FR,K&P} for Wirtinger presentations 
of knot quandles defined from
knot diagrams, which 
are 
similar to Wirtinger presentations of knot groups.
In this case, arcs of knot diagrams represent generators, and crossings
give relations of Wirtinger form.
Let 
$Q(K)$
represent such a quandle
for a knot diagram 
$K.$

\begin{sect}{\bf Example.\/}
{\rm 
If 
$K_0$ 
is the diagram illustrated
 in Fig.~\ref{ssvsv1}, then its quandle has presentation 
$$Q(K_0) =<x,y,z: x*y=z,y*z=x,y*x=z>.$$ 
The quandle
$<x,y,z: x*y=z,y*z=x,y*x=z>$ is isomorphic to 
the 3-element dihedral quandle $R_3$.
The abstract knot 
$K_0$ 
represents the 
$2$-cycle
$(x,y)+(y,z)-(y,x) \in Z^{\rm Q}_2( Q(K_0))$.
}
\end{sect}

\begin{sect}{\bf Proposition.} 
Every abstract  closed 
$1$-knot diagram 
$K$ 
represents a cycle   $[{K}] \in Z^{\rm Q}_2 (Q( K))$. 
\end{sect}

The abstract knot $K_0$ is a boundary since
the dihedral quandle, $R_3$, has trivial 
$2$-dimensional homology.

\begin{sect}{\bf Question.\/} 
Under what circumstances
 is the cycle $[K]$ a boundary? 
\end{sect}

\begin{sect}{\bf Discussion and Example.\/}  {\rm
On the other hand, we may define a shadow quandle for an abstract diagram.
We can show that the $3$-dimensional homology class 
that is represented by a shadow coloring of  Fig.~\ref{ssvsv1} is 
non-trivial as follows. The 5-element quandle $QS(5)$ (consisting 
of non-identity permutations in the symmetric group on 3 letters) 
colors the diagram non-trivially. Transpositions 
go on the arcs and $3$-cycles go on the $2$-dimensional 
regions. By a {\sc Mathematica} calculation, we have determined that the
$3$-cycle represented thereby is not a boundary.

}\end{sect}

\begin{figure}
\begin{center}
\mbox{
\epsfxsize=3in
\epsfbox{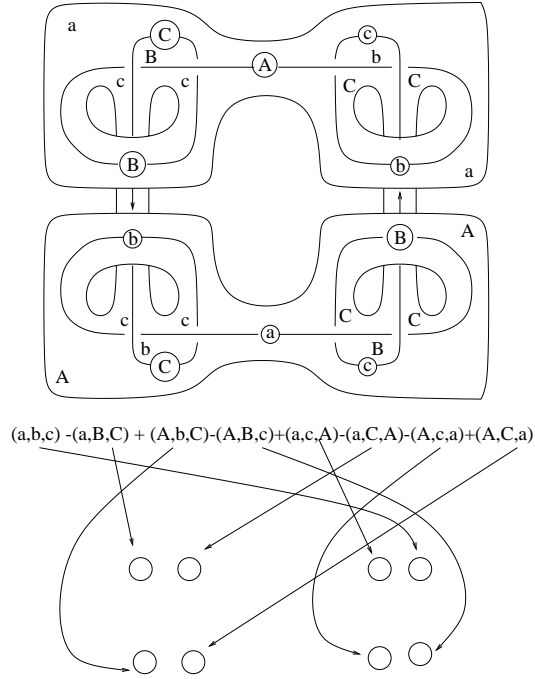}
}
\end{center}
\caption{A cycle of $QS(6)$ }
\label{s6cycle}
\end{figure}

\begin{figure}
\begin{center}
\mbox{
\epsfxsize=4in
\epsfbox{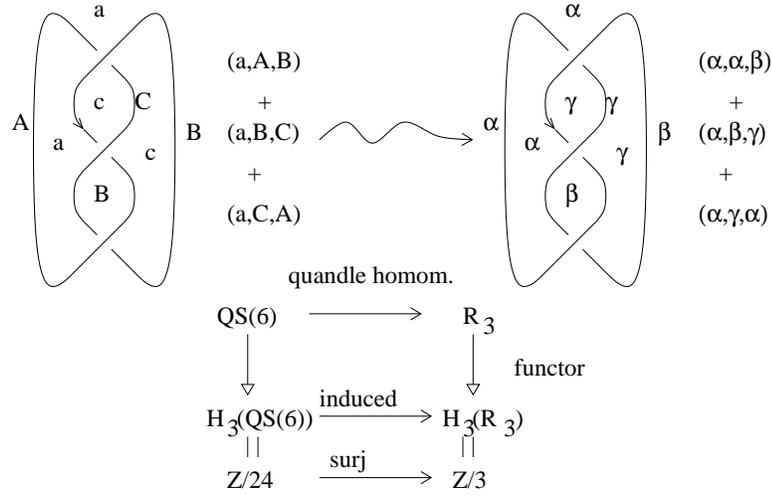}
}
\end{center}
\caption{A homology generator and its image under the induced map }
\label{genh3s6}
\end{figure}

\begin{sect}{\bf Example.\/} \label{qs62r3} 
{\rm The example depicted in Fig.~\ref{s6cycle} represents a 
$3$-cycle over the quandle $QS(6)$. The homomorphism 
$p:QS(6) \rightarrow R_3$  defined in \ref{quanxam}
induces a map on homology.
 The corresponding $R_3$ shadow 
colored  diagram represents a trivial $3$-cycle. We can show,
 by mean of a {\sc Mathematica} calculation, that 
the $QS(6)$ colored cycle is non-trivial.

The same program show that the colored diagram on the left of 
Fig.~\ref{genh3s6} represents a generator of $H^{\rm Q}_3(QS(6);{\bf Z}) = {\bf Z}_{24}$. The induced map $f_*$ is illustrated to be surjective.
} \end{sect}

\begin{sect} {\bf Remark.\/} {\rm
Let $K$ be a knot diagram on a compact oriented surface $F$.
Then the {\it fundamental 
shadow quandle} $SQ(K)$ is defined as follows. 
The generators correspond to 
over-arcs and connected components 
of $F \setminus \mbox{universe of} \ K$. 
The relations are defined for 
each crossing as ordinary 
fundamental quandles, and at each arc dividing 
regions. Specifically, if $a$ and 
$b$ are generators corresponding to
adjacent regions such that the 
normal points from the region colored $a$
to that colored $b$, and if the arc dividing these 
regions is colored by $c$, then we have the relation $b=a*c$.
This defines a presentation of a quandle, which is called the
fundamental shadow quandle of $K$. 
Two diagrams on $F$ that differ by Reidemeister moves on $F$
have 
isomorphic fundamental shadow quandles. 
The shadow colors are regarded as  quandle homomorphisms
from the fundamental shadow quandle to a 
quandle $X$.

} \end{sect}

\section{Boundary Homomorphisms $H^Q_{n+1} \rightarrow H^D_n$ }
\label{bdryhomsec}

In \cite{CJKS2}, from a split short exact sequence
\begin{eqnarray}
0 \to C_n^{\rm D}(X) \stackrel{i}{\to} C_n^{\rm R}(X) \stackrel{j}{\to}
C_n^{\rm Q}(X) \to 0,
\end{eqnarray}
the following homology long exact sequence 
\begin{eqnarray}
\cdots \stackrel{\partial_\ast}{\to} H_n^{\rm D}(X;G) \stackrel{i_\ast}{\to}
H_n^{\rm R}(X;G)
\stackrel{j_\ast}{\to} H_n^{\rm Q}(X;G)
\stackrel{\partial_\ast}{\to} H_{n-1}^{\rm D}(X;G) \to \cdots
\end{eqnarray}
was constructed. 
We give an application of 
Theorem~\ref{representation}
to boundary homomorphisms
of this exact sequence.

The endpoints  
of 
shadow 
colored arc diagrams are 
oriented by the orientation of the arcs. 
Each such endpoint represents a rack $2$-chain
$\pm(a,a)$ where $a$ is the color on the surrounding region.
The sum of these represents the image of the 
boundary map $\partial_* :H^{\rm Q}_3(X) \rightarrow H^{\rm D}_2(X)$
in the long exact  homology sequence.

Similarly, each oriented  
shadow 
colored
branch point represents a rack 
$3$-chain $\pm(a,b,b)$ where $a$ is the color
in one of the surrounding regions and $b$ 
is the color on the local surface. 
An
 oriented 
shadow 
colored 
hem $2$-crossing diagrams represent $3$-chains of the form
$(a,a,b)$ where $a$ is the color 
on one side of the upper sheet, 
and $b$ is the color on the upper sheet.
the image of the boundary map 
$\partial_* :H^{\rm Q}_4(X) 
\rightarrow H^{\rm D}_4(X)$ is represented as the sum
over all branch points and hem $2$-crossings of the chains 
these represent. 

In Theorems~\ref{parttrivtoo} and \ref{partialtriv}, 
we will use these descriptions 
 of the homomorphisms
$\partial_*$  and geometric techniques 
to show that these boundary maps 
are trivial.

\begin{figure}
\begin{center}
\mbox{
\epsfxsize=3in
\epsfbox{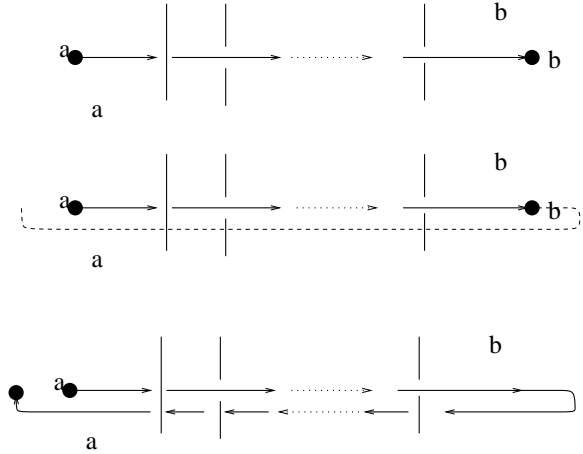}
}
\end{center}
\caption{The boundary homomorphism is trivial}
\label{trivquanbound}
\end{figure}

\begin{sect} {\bf Theorem \cite{CJKS2}.\/}  
\label{parttrivtoo}
 Let $X$ be a quandle.
The boundary homomorphism
 $ \partial_*  : H_3^Q(X) \rightarrow H_2^D(X)$
in the long exact sequence of quandle homology is trivial. 
\end{sect} 
{\it Proof.\/}
We give two proofs.

(1)
Let $\eta= \sum_{i=1}^k \epsilon_i (x_i, y_i, z_i) \in Z_3^Q(X)$. 
Then $\partial \eta = \sum_{j=1}^h 
\epsilon_j
(w_j, w_j)$
for some $w_j \in X$ for $j=1, \cdots, h$,
from the definition of quandle cycle groups.
Other terms coming out from each $\partial ( x_i, y_i, z_i)$
cancel out. 

Let $\tau_i$, $i=1, \cdots, k$, be the $3$-crossing diagrams
colored with the triple  $(x_i, y_i, z_i)$ for each $i$.
According to the cancelation of the terms  $\partial ( x_i, y_i, z_i)$,
paste together $\tau_i$. 
There are boundary crossings (boundaries of $\tau_i$)
that are colored by $(w_j, w_j)$.
As the double curves form an immersed 
$1$-manifold with boundary,
these boundary crossings are paired as
$\sum_{g=1}^{h/2}[ (a_g, a_g) - (b_g, b_g) ]$,
where the negative sign for $b_g$ represents that the orientation 
of the double curve points into 
 the corresponding crossing.
As one traces the double curve from the crossing with color
$ (a_g, a_g) $ to that with  $(b_g, b_g)$, 
the curve under-goes the middle or top sheets at some of $\tau_i$'s. 
The color changes accordingly, but we have that $a_g$ and $b_g$
belong to the same orbit (as $b_g = a_g * w_g$ for some word $w_g$ in $X$).

On the other hand, one computes that 
$\partial (a,a,b)=-(a,a)+(a*b, a*b)$, so that 
$(a_g, a_g) - (b_g, b_g) \in B_2^D(X)$, as desired.

(2) Represent the homology class by a cycle $\eta$ and choose
 a shadow colored 
abstract $1$-knot
diagram to represent $\eta$. Along each arc there are
 a collection of crossings.
These represent the summands of $\eta$. We push the
arc back along a parallel as indicated in Fig~\ref{trivquanbound}.
 When pushing backwards we make sure that the endpoint always passes under
 the crossings. Each  $2$-crossing that is introduced in the process is shadow colored of the form $(x,x,y)$.
 Thus the new crossings do not affect the 
quandle homology class represented by the diagram. The shadow colored endpoints can be pairwise canceled. Thus  $\partial_*[\eta] =0.$  This completes the proof. $\Box$

\begin{figure}
\begin{center}
\mbox{
\epsfxsize=6in
\epsfbox{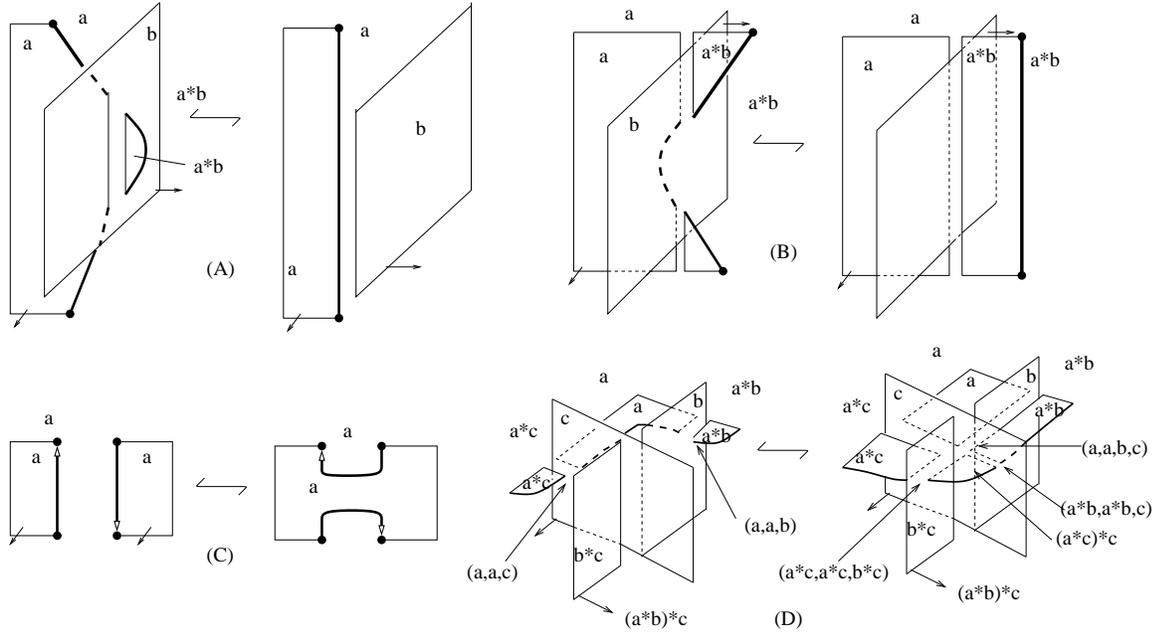}
}
\end{center}
\caption{ Moves of hems}
\label{funnystuff}
\end{figure}

\begin{sect}{\bf Theorem.} \label{partialtriv} 
 Let $X$ be a quandle.
The boundary homomorphism
$\partial_*: H^{\rm Q}_4 (X) \rightarrow H^{\rm D}_3 (X)$
in the long exact sequence of quandle homology is trivial.
\end{sect}
{\it Proof.} Consider a homology class 
$[\eta]\in H^{\rm Q}_4(X)$ and represent 
it by the cycle $\eta$. Construct a
shadow colored 
abstract 
$2$-knot 
diagram 
${\mathcal SD}_\eta
= [f : M \rightarrow N ] $  with 
hems to represent $\eta$. 
We may assume (by attaching handles 
if necessary) that the 
ambient  
3-manifold, 
$N$,
is closed. 
Moreover, the boundary of 
the surface $M$
consist 
of simple closed curves formed by the 
hems in the diagram.

The image
$\partial_*[\eta]$ is  represented on
${\mathcal SD}_\eta$ as the collection of 
shadow colored branch point
diagrams 
(terms of the form $(a,b,b)$) and
hemmed $2$-crossing diagrams 
(terms of the form
$(a,a,b)$). We first eliminate the terms of the form $(a,a,b)$.

Consider the collection, $B^1$ of hems. These form a closed $1$-manifold
in the $3$-manifold, $N$.  
The $1$-manifold $B^1$ is null homologous 
in 
$N$
since it is the boundary of the 
surface $M$ of ${\mathcal SD}_\eta$. 
Therefore $B^1$ bounds an embedded Seifert surface in $N$. The Seifert surface can be assumed to intersect the 
surface $M$ in general position. 
Now decompose the Seifert
surface into $1$- and $2$-handles 
where the $1$-handles are attached along a neighborhood of the hem. We will show that the hem can be eliminated by attaching the Seifert 
surface to the diagram ${\mathcal SD}_\eta$ 
along the hem. In the process, 
triple points will be introduced, 
but these will not affect the quandle homology class represented.

The core disk of a $1$-handle  of the Seifert surface 
intersects $M$ at a finite number of points. We push  a 
segment of 
the hem along this core disk  using the 
move
depicted  in (A) of 
Fig.~\ref{funnystuff} 
until segments are 
in a small ball neighborhood that contains no further sheet of the surface. 
The situation is depicted in (C) left 
 in Fig.~\ref{funnystuff}, 
and connect the segments as indicated in (C). 
 Having pushed the 
$B^1$ and joined up segments, 
we may assume that each component of the 
hem is unknotted and the collection is unlinked. 
Then an unknotting disk for a component 
of $B^1$ intersects the double point 
set of ${\mathcal SD}_\eta$ in a finite number of points. 
We push the hem
across these points introducing 
triple points, 
as indicated in Fig.~\ref{funnystuff} (D), 
that represent 
chains of the form $(a,a,b,c)$. 
(Recall that the hemmed sheet is always the under-sheet at each
hemmed $2$-crossing diagram).

But $[\eta \pm (a,a,b,c)] = 
[\eta] \in H^{\rm Q}_4(X)$.
The number of intersections with the 
resulting hem and $M$ is 
necessarily even since the hem is null homologous. 
We can push the components 
further, 
using the move depicted in Fig.~\ref{funnystuff} (B), 
until each bounds a disk that 
does not intersect the surface $M$.   
Now attach these  disks to $M$ to obtain a shadow colored 
diagram of a closed surface 
that represents $\eta$. In this way we may assume that
$\partial_*[\eta]$ 
consists entirely of  terms of the form
$(a,b,b)$.

The terms of the form $(a,b,b)$ are 
represented on the diagram by branch points. There are an even number of these since each arc of double points has 
two ends. Using the same technique 
as in \cite{CS:cancel},  
we can  cancel these in pairs.  
The process of moving the surface 
involves quandle colored Roseman moves.
The cancelation may introduce shadow colored 
triple points,  but  a consecutive pair of 
two of the four colors adjacent to the  
triple points introduced agree. Thus the triple points do not affect the quandle homology class that the diagrams represent.

In this way $\eta$ is represented by a shadow colored diagram of a closed surface without branch points. Thus $\partial_*[\eta]=0$. $\Box$

\noindent
{\bf Acknowledgements.} JSC is being supported by NSF grant DMS-9988107.
MS is being supported by NSF grant DMS-9988101. 
SK is being supported by a Fellowship from the Japan Society for the Promotion of Science. 
We have had productive conversations with Dan Silver about this paper.

\end{document}